\documentclass[preprint,3p,times]{elsarticle}
\usepackage[toc,page,title,titletoc,header]{appendix}
\usepackage{amsmath}
\usepackage{amssymb}
\usepackage{amsthm}
\usepackage{graphicx}
\usepackage{epsfig}
\usepackage{tikz,pgfplots,tikz-3dplot}
\usepackage{epstopdf}
\usetikzlibrary{fit}
\biboptions{sort&compress}

\newtheorem{thm}{Theorem}[section]
\newtheorem{lemma}[thm]{Lemma}

\newcommand{\vnorm}[1]{\left\lVert#1\right\rVert}
\newcommand{\bRplus}{{\mathbb R}_{>0}}
\newcommand{\bRgeq}{{\mathbb R}_{\geq 0}}
\newcommand{\RZ}{{\mathbb R} \slash {\mathbb Z}}
\newcommand{\bR}{{\mathbb R}}

\newcommand{\dA}{\;{\rm{d} A}}

\newcommand{\vol}{\operatorname{vol}}

\newcommand{\mbI}{\mathbb{I}}
\newcommand{\rd}{\;{\rm d}}

\newcommand{\drho}{\;{\rm d}\rho}

\newcommand{\Vhpartial}{\underline{V}^h_\partial}
\newcommand{\Vhpartialzero}{\underline{V}^h_{\partial_0}}
\newcommand{\Vpartial}{\underline{V}_\partial}
\newcommand{\Vpartialzero}{\underline{V}_{\partial_0}}

\newcommand{\id}{\rm id}

\newcommand{\dd}[1]{\frac{\rm d}{{\rm d}#1}}
\newcommand{\ddt}{\dd{t}}

\newcommand{\nn}{\nonumber}
\newcommand{\ek}{e}
\newcommand{\ttau}{\Delta t}
\newcommand{\sliprho}
{\widehat\varrho_{_{\partial\mathcal{S}}}}
\newcommand{\normal}{{\rm n}}
\def\epsilon{\varepsilon}

\newcommand{\pol}{\vec}  

\begin{document}
\begin{frontmatter}
\title{
Volume-preserving parametric finite element methods for 
axisymmetric geometric evolution equations
}

\author[1]{Weizhu Bao}
\address[1]{Department of Mathematics, National University of Singapore, 119076, Singapore}
\ead{matbaowz@nus.edu.sg}
\author[2]{Harald Garcke}
\address[2]{Fakult{\"a}t f{\"u}r Mathematik, Universit{\"a}t Regensburg, 
93040 Regensburg, Germany}
\ead{harald.garcke@ur.de}
\author[3]{Robert N\"urnberg}
\address[3]{Dipartimento di Mathematica, Universit\`a di Trento,
38123 Trento, Italy}
\ead{robert.nurnberg@unitn.it}
\author[2]{Quan Zhao}
\ead{quan.zhao@ur.de}

\begin{abstract}
We propose and analyze volume-preserving parametric finite element methods for 
surface diffusion, conserved mean curvature flow and an intermediate evolution law in an axisymmetric setting. The weak formulations are presented in terms of the generating curves of the axisymmetric surfaces. The proposed numerical methods are based on piecewise linear parametric finite elements. The constructed fully practical schemes satisfy the conservation of the enclosed volume. 
In addition, we prove the unconditional stability and consider the distribution of vertices for the discretized schemes. The introduced methods are implicit and the resulting nonlinear systems of equations can be solved very efficiently and accurately via the Newton's iterative method. Numerical results are presented to show the accuracy and efficiency of the introduced schemes for computing the considered axisymmetric geometric flows.
\end{abstract} 

\begin{keyword} Surface diffusion flow, conserved mean curvature flow,
parametric finite element method, axisymmetry, volume conservation, 
unconditional stability, Newton's method
\end{keyword}

\end{frontmatter}

\setcounter{equation}{0}
\section{Introduction} \label{sec:intro}

The motion of interfaces driven by a law for the normal velocity which involves curvature quantities plays an important role in applied mathematics and materials science. One of the most prominent examples is the surface diffusion flow, which was first proposed by Mullins to describe the evolution of microstructure in polycrystalline materials \cite{Mullins57}. In this evolution law the normal velocity of the interface is given by the surface Laplacian of the mean curvature, and the resulting differential equations are parabolic and of fourth order \cite{Mullins57,Cahn94}. Let $\left\{\mathcal{S}(t)\right\}_{t \geq 0} \subset \bR^3$ be a family of smooth, oriented hypersurfaces, which for now we assume to be closed. The motion by surface diffusion flow is then given by 
\begin{equation} \label{eq:sdS}
\mathcal{V}_{n} =- \Delta_{_\mathcal{S}}\,\mathcal{H} \qquad\text{on }
\mathcal{S}(t)\,,
\end{equation}
where $\mathcal{V}_{n}$ denotes the normal velocity of 
$\mathcal{S}(t)$ in the direction of the normal $\pol\normal_{_\mathcal{S}}$, $\Delta_{_\mathcal{S}}$ is the Laplace--Beltrami operator on
$\mathcal{S}(t)$, and $\mathcal{H}$ is the mean curvature of $\mathcal{S}(t)$, given by \cite{Deckelnick05}
\begin{equation*} 
\Delta_{_\mathcal{S}}\,\pol\id = \mathcal{H} \,
\pol\normal_{_\mathcal{S}}  \qquad \text{ on } \mathcal{S}(t)\,,
\end{equation*} 
with $\pol\id$ being the identity function in $\mathbb{R}^3$. A geometric flow  which combines surface diffusion and surface attachment limited kinetics, which we call the intermediate evolution flow, is given by
\begin{equation} \label{eq:imS}
\mathcal{V}_{n} =-\Delta_{_\mathcal{S}}\, \mathcal{Y}\,, \qquad \left(
-\frac1\xi\, \Delta_{_\mathcal{S}}+ \frac1\alpha \right) \mathcal{Y} =\mathcal{H}
\quad\text{on }\ \mathcal{S}(t) \,,
\end{equation}
where $\alpha,\xi \in \bRplus$ are given positive parameters. In the limit 
of fast attachment kinetics $\xi\to \infty$ and $\alpha=1$, formal arguments
suggest, see \cite{Cahn94,Taylor94}, that one recovers the surface diffusion 
flow (\ref{eq:sdS}).
In the limit of fast surface diffusion $\alpha\to\infty$
and $\xi=1$, one expects, see \cite{Cahn94,Taylor94}, to recover the 
conserved mean curvature flow,
\begin{equation} \label{eq:cmcfS}
\mathcal{V}_{n} = \mathcal{H}
- \frac{\int_{\mathcal{S}(t)} \mathcal{H} \dA}{\int_{\mathcal{S}(t)} 1 \dA}
\quad\text{on }\ \mathcal{S}(t)\,,
\end{equation}
where dA is the surface measure. We note that all three flows can be
interpreted as volume conserving gradient flows for the surface area, since
they satisfy
\begin{subequations}\label{eq:dtareavol}
\begin{equation} \label{eq:dtarea}
\ddt \int_{\mathcal{S}(t)} 1\dA = -
\int_{\mathcal{S}(t)} \mathcal{H}\,\mathcal{V}_n \dA \leq 0
\end{equation}
and
\begin{equation} \label{eq:dtvol}
\ddt \vol(\Omega(t)) = \int_{\mathcal{S}(t)} \mathcal{V}_n \dA = 0,
\end{equation}
\end{subequations}
if $\Omega(t)$ denotes the interior of $\mathcal{S}(t)$ and 
$\pol\normal_{_\mathcal{S}}$ is the outer normal to $\Omega(t)$ on 
$\partial\Omega(t)=\mathcal{S}(t)$.
Existence results for the flows \eqref{eq:sdS}, \eqref{eq:imS} and 
\eqref{eq:cmcfS} can be found in 
\cite{Gage86,Elliott1994,EscherMS98,Huisken1987volume}.
In \cite{Gage86,Huisken1987volume} it is shown that convexity is conserved for
conserved mean curvature flow. However, for surface diffusion and the 
intermediate evolution flow 
convexity in general is not conserved, see \cite{GigaI99,Ito02,EscherI05}.
It is rigorously shown in \cite{EscherGI02} that for $\alpha\to\infty$ and
$\xi=1$ the intermediate flow \eqref{eq:imS} converges to the conserved mean
curvature flow. So far it is an open problem to rigorously show that for
$\xi\to\infty$ and $\alpha=1$ the intermediate flow converges to motion by
surface diffusion. 
For a more detailed discussion, we refer the reader to Refs.~\cite{Taylor94,Cahn94,Elliott1994,Barrett08JCP,Barrett19} and the references therein.  

In this paper, we will investigate the numerical approximation of the three
flows \eqref{eq:sdS}, \eqref{eq:imS} and \eqref{eq:cmcfS}, paying particular
attention to the volume preserving aspect of these flows. 
Numerical approximations for curvature driven flows have been studied intensively during the last 30 years, e.g., \cite{Dziuk90,ColemanFM96,Dziuk2002,Bansch04, Bansch05,Barrett07, Barrett07b,Barrett08,Barrett2011, Kovacs2019, Kovacs2020,Li2020energy, Jiang2021, Zhao2020p,Zhao2021}, and we refer to the review articles \cite{Deckelnick05, Barrett20} for detailed discussions. Earlier front tracking methods often require mesh regularization\slash smoothing algorithms to prevent mesh degenerations during the simulations \cite{Dziuk2002,Bansch05, Wang15}. Based on ideas of Dziuk \cite{Dziuk2002}, Barrett, Garcke, and N\"urnberg (BGN) first presented variational approximations based on weak formulations that allow tangential degrees of freedom so that the mesh quality can be improved significantly for the introduced parametric finite element methods \cite{Barrett07b, Barrett07, Barrett08JCP, Barrett20}. In particular, for the evolution of curves, a semi-discrete version of their discretization method leads to the equidistribution of mesh points \cite{Barrett07}, while in the case of surfaces their semi-discrete schemes lead to so-called conformal polyhedral surfaces \cite{Barrett08JCP}. In both cases the observed tangential motion of vertices in practice helps to avoid the complex and undesired re-meshing steps. Moreover, the fully discrete BGN schemes usually enjoy unconditional stability, which mimics the surface area dissipation law in \eqref{eq:dtarea} for the flows we consider in this paper. However, these fully discrete schemes in general fail to exactly conserve the volume enclosed by the discrete surfaces. The observed volume loss for these discretized schemes can be significant and non-negligible under certain circumstances. This could pose a great challenge on the accuracy of the numerical solutions, and thus numerical schemes that can conserve the enclosed volume are desirable and necessary.  

In a separate development, very recently Bao and Zhao, based on the original
ideas of BGN, have introduced a numerical method for surface diffusion flow
that satisfies an exact volume conservation property on the fully discrete 
level \cite{Zhao2021}. Here the crucial idea is to use a suitably
time-integrated discrete normal vector, an idea that was first explored in the
context of surface diffusion for curves in the plane, and in the absence of any
tangential motion, in \cite{Jiang2021}. It is the aim of this paper to
transfer the novel approach in \cite{Jiang2021,Zhao2021} to the approximation
of the volume preserving flows \eqref{eq:sdS}, \eqref{eq:imS} and 
\eqref{eq:cmcfS} in the axisymmetric setting.


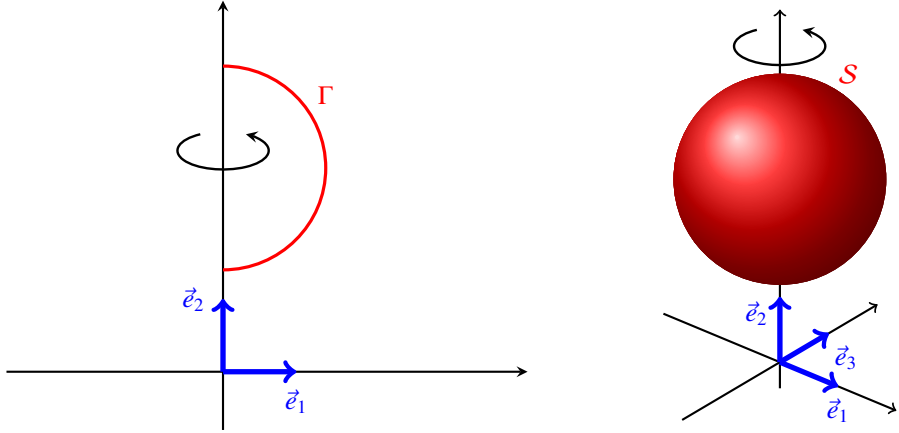
\begin{figure}
\center
\newcommand{\AxisRotator}[1][rotate=0]{%
    \tikz [x=0.25cm,y=0.60cm,line width=.2ex,-stealth,#1] \draw (0,0) arc (-150:150:1 and 1);%
}
\begin{tikzpicture}[every plot/.append style={very thick}, scale = 1]
\begin{axis}[axis equal,axis line style=thick,axis lines=center, xtick style ={draw=none}, 
ytick style ={draw=none}, xticklabels = {}, 
yticklabels = {}, 
xmin=-0.2, xmax = 0.8, ymin = -0.4, ymax = 2.55]
after end axis/.code={  
   \node at (axis cs:0.0,1.5) {\AxisRotator[rotate=-90]};
   \draw[blue,->,line width=2pt] (axis cs:0,0) -- (axis cs:0.5,0);
   \draw[blue,->,line width=2pt] (axis cs:0,0) -- (axis cs:0,0.5);
   \node[blue] at (axis cs:0.5,-0.2){$\pol\ek_1$};
   \node[blue] at (axis cs:-0.2,0.5){$\pol\ek_2$};
   \draw[red,very thick] (axis cs: 0,0.7) arc[radius = 70, start angle= -90, end angle= 90];
   \node[red] at (axis cs:0.7,1.9){$\Gamma$};
}
\end{axis}
\end{tikzpicture} \qquad \qquad
\tdplotsetmaincoords{120}{50}
\begin{tikzpicture}[scale=2, tdplot_main_coords,axis/.style={->},thick]
\draw[axis] (-1, 0, 0) -- (1, 0, 0);
\draw[axis] (0, -1, 0) -- (0, 1, 0);
\draw[axis] (0, 0, -0.2) -- (0, 0, 2.7);
\draw[blue,->,line width=2pt] (0,0,0) -- (0,0.5,0) node [below] {$\pol\ek_1$};
\draw[blue,->,line width=2pt] (0,0,0) -- (0,0.0,0.5);
\draw[blue,->,line width=2pt] (0,0,0) -- (0.5,0.0,0);
\node[blue] at (0.2,0.4,0.1){$\pol\ek_3$};
\node[blue] at (0,-0.2,0.3){$\pol\ek_2$};
\node[red] at (0.7,0,1.9){$\mathcal{S}$};
\node at (0.0,0.0,2.4) {\AxisRotator[rotate=-90]};

\tdplottransformmainscreen{0}{0}{1.4}
\shade[tdplot_screen_coords, ball color = red] (\tdplotresx,\tdplotresy) circle (0.7);
\end{tikzpicture}
\caption{Sketch of $\Gamma$ and $\mathcal{S}$, as well as 
the unit vectors $\pol\ek_1$, $\pol\ek_2$ and $\pol\ek_3$.}
\label{fig:sketch}
\end{figure}

In fact, in many situations the three-dimensional evolving surface may appear rotationally symmetric so that the studied geometric flows can be reduced to one-dimensional problems \cite{Bernoff1998, Deckelnick03, Jiang19b, Zhao19, Barrett19, Barrett2019variational, Barrett2021stable}. In addition, numerical schemes for the axisymmetric problem often avoid issues related to mesh distortions. Hence in this paper, we consider the flows \eqref{eq:sdS}, \eqref{eq:imS} and \eqref{eq:cmcfS} in the case that $\mathcal{S}(t)$ is an axisymmetric surface that is rotationally symmetric with respect to the $x_2$--axis, as shown in Fig.~\ref{fig:sketch}. In particular, we aim to propose volume-preserving parametric finite element methods by combining the ideas from \cite{Jiang2021, Zhao2021} with the recently introduced axisymmetric variational discretizations for the three flows from \cite{Barrett19,Barrett2019variational}. By using carefully weighted discrete normals, we are able to derive schemes that combine the volume preserving properties from \cite{Zhao2021} with the unconditional stability and good mesh properties from \cite{Barrett19,Barrett2019variational}. The obtained schemes are implicit, but the resulting systems of nonlinear equations can be accurately and efficiently solved via Newton's method. 

So far, for simplicity, we have assumed that the surfaces $\mathcal{S}(t)$ have
no boundary. But for many physical applications it is of interest to also
consider the case of surfaces with boundary, for example in the study of 
capillary surfaces \cite{Finn86,Barrett10} and in the study of the
solid-state dewetting of thin films on a substrate \cite{Zhao19,JiangZB20}. 
Hence in this paper we allow $\mathcal{S}(t)$ to be with or without boundary.
We assume that $\mathcal{S}(t)$ is made up of a single connected component, 
which due to the axisymmetric nature of $\mathcal{S}(t)$ implies that its
boundary, unless it is empty, consists of either one or two circles
that each lie within a hyperplane that is parallel to the $x_1$-$x_3$--plane. 
On writing the boundary as a disjoint union $\partial\mathcal{S}(t) = 
\partial_{C_1}\mathcal{S}(t) \cup \partial_{C_2} \mathcal{S}(t)$, where the
sets $\partial_{C_i}\mathcal{S}(t)$ are allowed to be empty, physically
meaningful boundary conditions for the flows we consider can then be formulated
as follows. In the simplest case, the boundary part 
$\partial_{C_i}\mathcal{S}(t)$ is fixed:
\begin{subequations} \label{eq:bc}
\begin{equation} \label{eq:Dbc}
\partial_{C_i}\mathcal{S}(t) = \partial_{C_i}\mathcal{S}(0) \qquad t \geq 0\,.
\end{equation}
Alternatively, we allow the boundary part $\partial_{C_i}\mathcal{S}(t)$ 
to either move with a velocity that is perpendicular to the
$x_1$-$x_3$--plane, or to expand or shrink in a hyperplane that is parallel
to the $x_1$-$x_3$--plane:
\begin{align} 
\partial_{C_i}\mathcal{S}(t) \subset 
\{ \pol z \in \bR^3 : (\pol z \cdot \pol\ek_1)^2 + (\pol z \cdot \pol\ek_3)^2
= r_i^2\}\,,\quad r_i > 0\,,\qquad t \geq 0\,, \label{eq:free_slip1} \\
\partial_{C_i}\mathcal{S}(t) \subset 
\{ \pol z \in \bR^3 : \pol z \cdot \pol\ek_2 = h_i\}\,,\quad h_i \in\bR
\,,\qquad t \geq 0\,. \label{eq:free_slip2}
\end{align}
\end{subequations}
The energy decay \eqref{eq:dtarea} then still holds for the three evolution
equations \eqref{eq:sdS}, \eqref{eq:imS} and \eqref{eq:cmcfS}, if in the case
of the moving boundaries \eqref{eq:free_slip1}, \eqref{eq:free_slip2} one
postulates a $90^\circ$ contact angle with the external walls/substrates. 
Moreover, the two flows \eqref{eq:sdS} and \eqref{eq:imS} need 
additional conditions on $\partial\mathcal{S}(t)$ for closure, e.g.\ no-flux
conditions for $\mathcal{H}$ or $\mathcal{Y}$. 
Similarly, the volume conservation 
\eqref{eq:dtvol} also still holds in the presence of the boundary conditions
\eqref{eq:bc}, on defining the ``interior'' $\Omega(t)$ 
of $\mathcal{S}(t)$ as a suitable
finite domain bounded by $\mathcal{S}(t)$ and parts of the external walls and
substrates.
Physically the $90^\circ$ contact angle condition
arises for a neutral external boundary. In situations where the interior and
exterior of $\mathcal{S}(t)$ have different contact energy densities with the
external boundary, their
difference multiplied with the area of the contact patch of $\Omega(t)$ with
the external boundary enters the free energy, and a generalized version
of \eqref{eq:dtarea} then still holds on prescribing the correct boundary
contact angles.
The precise conditions are described in e.g.\ \cite{Barrett10}, and we will
present the details in the axisymmetric setting in the next section.

The rest of the paper is organized as follows. In Section~\ref{sec:wf}, we 
review the strong and weak formulations for the three geometric evolution equations that we consider in this paper, at first only for the case of a neutral external boundary. Based on the discretizations of these weak formulations, several volume-preserving parametric finite element methods are proposed in Section~\ref{sec:fd}. For these schemes the properties of volume conservation, stability and vertex distribution are also analyzed in Section~\ref{sec:fd}. At the end of the section, we briefly discuss the extension of the presented schemes to the case of a non-neutral external boundary. Subsequently, several numerical results are presented in Section~\ref{sec:nr}. Finally the paper is concluded in Section~\ref{sec:con}. 
 
\setcounter{equation}{0} 
\section{Mathematical formulations}
\label{sec:wf}

In this section we briefly review strong and weak formulations for the three
evolution laws \eqref{eq:sdS}, \eqref{eq:imS} and \eqref{eq:cmcfS} in the
axisymmetric setting. We refer the reader to 
\cite{Barrett19,Barrett2019variational} for more details and the precise
derivations.

Let $\RZ$ be the periodic unit interval, and set
\[
\mbI = \RZ\,,\text{ with } \partial \mbI = \emptyset\,,\qquad \text{or}\quad
\mbI = (0,1)\,, \text{ with } \partial \mbI = \{0,1\}\,.
\]
We consider the axisymmetric situation, where 
$\pol x(\cdot,t)=(r(\cdot,t),~z(\cdot,t))^T : \overline {\mbI} \to 
\bRgeq \times \bR$ is a parameterization of $\Gamma(t)$. 
Throughout $\Gamma(t)$ represents the generating curve of a
surface $\mathcal{S}(t)$ 
that is axisymmetric with respect to the $x_2$--axis, as shown in Fig.~\ref{fig:sketch}.  In particular, an induced parameterization of $\mathcal{S}(t)$ is 
given by
\begin{equation*}  
(\rho,~\theta,~t)\mapsto 
\bigl(r(\rho,~t)\,\cos\theta, ~z(\rho,t), ~r(\rho,t)\,\sin\theta\bigr)^T, 
\quad \rho\in\mbI, \quad\theta \in [0,2\pi].
\end{equation*}
We allow $\Gamma(t)$ to be either a closed curve, parameterized over
$\RZ$, which corresponds to $\mathcal{S}(t)$ being a genus-1 surface
without boundary.
Or $\Gamma(t)$ may be an open curve, parameterized over $[0,1]$. 
If both endpoints of $\Gamma(t)$ lie on the $x_2$--axis, then the surface
$\mathcal{S}(t)$ is a genus-0 surface without boundary. In all other cases
the boundary $\mathcal{S}(t)$ is made up of one or two circles that each
lie in a hyperplane that is parallel to the 
$x_1$-$x_3$--plane $\bR \times \{0\} \times \bR$.
Overall, we assume that the parameterization $\pol x$ satisfies the following
conditions, for all $t \in [0,T]$:
\[
\pol x(\rho,t) \cdot\pol\ek_1 > 0 \quad 
\forall\ \rho \in \overline {\mbI}\setminus \partial_0 \mbI\,,
\]
as well as
\begin{subequations} \label{eq:bcs}
\begin{align} 
\pol x(\rho,t) \cdot\pol\ek_1 &= 0 \quad 
\forall\ \rho \in \partial_0 \mbI\,,\label{eq:axibc} \\
\pol x_t(\rho,t) \cdot\pol\ek_i &= 0 \quad 
\forall\ \rho \in \partial_i\mbI \,, \ i =1,2\,,  \label{eq:freeslipbc} \\
\pol x_t(\rho,t) &= \pol 0 \quad 
\forall\ \rho \in \partial_{_D} \mbI \,, \label{eq:noslipbc}
\end{align}
\end{subequations}
where $\partial_{_D}\mbI \cup \bigcup_{i=0}^2 \partial_i\mbI = \partial\mbI$ 
is a disjoint partitioning of $\partial\mbI$. 
Here $\partial_{_D}\mbI \cup \bigcup_{i=1}^2 \partial_i\mbI$ 
denotes the subset of boundary points of $\mbI$ that model components of the
boundary of $\mathcal{S}(t)$, with $\partial_D\mbI$ corresponding
to fixed boundary circles as in \eqref{eq:Dbc}. 
Moreover, endpoints in $\partial_1\mbI$ correspond to boundary circles 
that can move freely up and down along the boundary of an infinite cylinder
that is aligned with the axis of rotation, \eqref{eq:free_slip1}, while 
$\partial_2\mbI$ corresponds to boundary circles
that can freely expand and shrink within a hyperplane that is parallel to
the $x_1$-$x_3$--plane, \eqref{eq:free_slip2}.
For a visualization of these different types of boundary nodes, the reader may 
refer to Table~1 in \cite{Barrett19}.

On assuming that $|\pol x_\rho(\cdot,t)| > 0$ in $\mbI$,
we introduce the arc length $s$ of the curve $\Gamma(t)$, i.e.\ $\partial_s =
|\pol{x}_\rho|^{-1}\,\partial_\rho$, and define the unit tangent and unit
normal to the curve $\Gamma(t)$ via
$\pol\tau = \pol x_s$ and $\pol\nu = -\pol\tau^\perp$,
where $(\cdot)^\perp$ denotes a clockwise rotation by $\frac{\pi}{2}$. 
We assume from now on that $\pol\nu(t)$ is such that the induced normal 
on $\mathcal{S}(t)$ is the {\em outer} normal to the ``interior'' $\Omega(t)$
enclosed by $\mathcal{S}(t)$ and possibly parts of the external boundary.
We also introduce the curvature $\kappa$ of the curve $\Gamma(t)$ via
$\kappa\,\pol\nu = \pol x_{ss}$, as well as the mean curvature of 
$\mathcal{S}(t)$ via $\varkappa = \mathcal{H} \circ \binom{\pol x}{0}$. 
Then it holds that
\begin{equation} \label{eq:kappaS}
\varkappa = \kappa
- \frac{\pol\nu\cdot\pol\ek_1}{\pol x\cdot\pol\ek_1}\,.
\end{equation}

For the geometric flows we consider in this paper, the surface area and the
enclosed volume play an important role.
If we denote by $A(\pol x(t))$ and $M(\pol x(t))$ the surface area and the total enclosed volume of the axisymmetric hypersurface $\mathcal{S}(t)$, respectively, then we have 
\begin{subequations}
\begin{align}
A(\pol x(t))&= \int_{\mathcal{S}(t)} 1\dA = 
2\pi\,\int_{\mbI} \pol x(\rho,t)\cdot\pol\ek_1\,|\pol x_\rho(\rho,t)|\drho,\label{eq:TE}\\
M(\pol x(t))&=\vol(\Omega(t))\,.\label{eq:defM}
\end{align}
\end{subequations}
Direct calculation, and recalling \eqref{eq:dtvol}, then yields that
\begin{subequations}
\begin{align} 
\ddt A(\pol x(t)) & = 2\pi\,\int_\mbI \left[\pol x_t\cdot\pol\ek_1
+ \pol x\cdot\pol\ek_1\,\frac{(\pol x_t)_\rho\cdot\pol x_\rho}{|\pol x_\rho|^2}
\right] |\pol x_\rho| \drho\,, \label{eq:dta} \\
\ddt M(\pol x(t)) & = 2\pi\int_{\mbI}(\pol x\cdot\pol\ek_1)\,\pol x_t\cdot\pol\nu \,|\pol x_\rho|\drho\,.\label{eq:dtm}
\end{align}
\end{subequations}

\subsection{Surface diffusion flow} \label{sec:sdf}

In the axisymmetric setting, the surface diffusion flow in \eqref{eq:sdS}
can be written in terms of an evolution equation for the generating curves $\Gamma(t)=\pol x(\cdot,~t)$ as follows:
\begin{equation} \label{eq:xtsd}
(\pol x\cdot\pol\ek_1)\,\pol x_t\cdot\pol\nu = 
- \left[\pol x\cdot\pol\ek_1\,\varkappa_s \right]_{s}
\quad\text{in }\ \mbI \times (0,T]\,,
\end{equation}
together with the boundary conditions \eqref{eq:bcs}
for $t\in(0,T]$, as well as
\begin{subequations} \label{eq:sdbdn}
\begin{alignat}{2}
\pol x_\rho\cdot\pol\ek_2&=0&&\quad\text{on }\partial_0\mbI \times(0,T],
\label{eq:bc4} \\
\pol x_\rho\cdot\pol\ek_{3-i}&=0&&\quad\text{on }\partial_i\mbI \times(0,T],\
i=1,2\,,\label{eq:bc4a} \\
\varkappa_\rho&=0&&\quad\text{on }\partial\mbI \times(0,T].\label{eq:bc5}
\end{alignat}
\end{subequations}
 
In order to introduce the weak formulations, we define the function spaces
\begin{align*}
\Vpartialzero&=\Bigl\{\pol\eta \in [H^1(\mbI)]^2 :\; \pol\eta(\rho)\cdot\pol\ek_1 = 0
\;\; \forall\ \rho \in \partial_0 \mbI\Bigr\},\\
\Vpartial &= \Bigl\{ \pol\eta\in\Vpartialzero :\, \pol\eta(\rho)\cdot\pol\ek_i = 0
\; \forall\ \rho \in \partial_i\mbI\,,\ i=1,2; \,
\pol\eta(\rho) = \pol 0 \;\; \forall\ \rho \in \partial_{_D}\mbI\Bigr\}.
\end{align*}
Then multiplying \eqref{eq:xtsd} by $\chi\,|\pol x_\rho|$ for
a test function $\chi\in H^1(\mbI)$, integrating over $\mbI$ 
and integrating by parts, on noting the boundary condition \eqref{eq:bc5}, 
yields
\[
\int_{\mbI}  (\pol x\cdot\pol\ek_1)\,\pol x_t\cdot\pol\nu\,\chi\,|\pol x_\rho|\drho
- \int_{\mbI}  \pol x\cdot\pol\ek_1\,
\varkappa_\rho\, \chi_\rho\,|\pol x_\rho|^{-1} \drho=0
\qquad \forall\ \chi\in H^1(\mbI).
\]
Similarly, it was shown in \cite{Barrett19} that multiplying \eqref{eq:kappaS} 
by $(\pol x\cdot\pol\ek_1)\pol\nu\cdot\pol\eta\,|\pol x_\rho|$ for a test function 
$\pol\eta\in \Vpartial$ leads to
\begin{equation} \label{eq:weakb}
 \int_{\mbI} (\pol x\cdot\pol\ek_1)\,
\varkappa\,\pol\nu\cdot\pol\eta\, |\pol x_\rho| \drho
+\int_{\mbI} \left[\pol\eta \cdot\pol\ek_1
+ \pol x\cdot\pol\ek_1\,\frac{\pol x_\rho\cdot\pol\eta_\rho}{|\pol x_\rho|^2}
\right] |\pol x_\rho| \drho 
= 0\qquad \forall\ \pol\eta\in \Vpartial\,.
\end{equation}
We mention that alternatively \eqref{eq:weakb} can also be directly deduced
from \eqref{eq:dtarea} and \eqref{eq:dta}. 

On defining $\langle\cdot,\cdot\rangle$ as the $L^2$-inner product over $\mbI$, we hence consider the following weak formulation for \eqref{eq:xtsd}.

\noindent
Let $\pol x(0) \in \Vpartialzero$. For $t \in (0,T]$
we find $\pol x(t) \in [H^1(\mbI)]^2$, with $\pol x_t(t) \in \Vpartial$, 
and $\varkappa(t)\in H^1(\mbI)$ such that
\begin{subequations} \label{eq:sdstabab}
\begin{alignat}{2}
 \Bigl\langle  (\pol x\cdot\pol\ek_1)\,\pol x_t\cdot\pol\nu,~\chi\,|\pol x_\rho|\Bigr\rangle
- \Big\langle \pol x\cdot\pol\ek_1\,
\varkappa_\rho, ~\chi_\rho\,|\pol x_\rho|^{-1} \Big\rangle&=0 &&
\qquad \forall\ \chi \in H^1(\mbI)\,, \label{eq:sdstaba} \\
\Big\langle (\pol x\cdot\pol\ek_1)\,
\varkappa\,\pol\nu,~\pol\eta\, |\pol x_\rho| \Big\rangle + \Big\langle\pol\ek_1,~\pol\eta\,|\pol x_\rho| \Big\rangle 
+\Big\langle
(\pol x\cdot\pol\ek_1)\,\pol x_\rho,~\pol\eta_\rho|\pol x_\rho|^{-1}\Big\rangle &=0 && \qquad \forall\ \pol\eta \in \Vpartial\,.
\label{eq:sdstabb}
\end{alignat}
\end{subequations}
Choosing $\chi = 1$ in \eqref{eq:sdstaba} and recalling \eqref{eq:dtm} yields 
that \eqref{eq:dtvol} is satisfied, while choosing $\chi=\varkappa$ in \eqref{eq:sdstaba} and $\pol\eta = \pol x_t$ in \eqref{eq:sdstabb}, on noting \eqref{eq:dta}, yields that
\begin{equation} \label{eq:dtasd}
\ddt A(\pol x(t)) = - 2\pi
\left\langle\pol x\cdot\pol\ek_1\,|\varkappa_\rho|^2,~|\pol x_\rho|^{-1}\right\rangle \leq 0\,,
\end{equation}
which shows that \eqref{eq:dtarea} holds. This implies the volume conservation and energy dissipation within the weak formulation. 

We now consider an alternative weak formulation, which treats the curvature $\kappa$ of the curve $\Gamma(t)$ as an unknown, and which will lead to an equidistribution property on the discrete level.
To this end, we let 
$\lambda=\frac{\pol\nu\cdot\pol\ek_1}{\pol x\cdot\pol\ek_1}$ and recall
that $\varkappa = \kappa - \lambda$.

\noindent
Let $\pol x(0) \in \Vpartialzero$. For $t \in (0,T]$
we find $\pol x(t) \in [H^1(\mbI)]^2$, with $\pol x_t(t) \in \Vpartial$, 
and $\kappa(t)\in H^1(\mbI)$ such that
\begin{subequations} \label{eq:sdweak}
\begin{alignat}{2}
\Big\langle(\pol x\cdot\pol\ek_1)\,\pol x_t\cdot\pol\nu,~\chi\,|\pol x_\rho|\Big\rangle - \left\langle  \pol x\cdot\pol\ek_1\, 
\left(\kappa - \lambda
\right)_\rho,~ \chi_\rho\,|\pol x_\rho|^{-1}\right\rangle &=0 &&
\qquad \forall\ \chi \in H^1(\mbI)\,, \label{eq:sdweaka} \\
\Big\langle \kappa\,\pol\nu,~\pol\eta\, |\pol x_\rho| \Big\rangle
+ \Big\langle \pol x_\rho,~\pol\eta_\rho\, |\pol x_\rho|^{-1} \Big\rangle
& = 0&& \qquad \forall\ \pol\eta \in \Vpartial\,,
\label{eq:sdweakb}
\end{alignat}
\end{subequations}
where \eqref{eq:sdweakb} is derived based on the curvature formulation $\kappa\,\pol\nu = \pol x_{ss}$. As before, choosing $\chi = 1$ in \eqref{eq:sdweaka} and noting \eqref{eq:dtm} yields that \eqref{eq:dtvol} is satisfied. Moreover, on choosing $\chi = \kappa -\lambda$ in \eqref{eq:sdweaka} and $\pol\eta = (\pol x\cdot\pol\ek_1)\,\pol x_t$ in \eqref{eq:sdweakb}, it is not difficult to prove that \eqref{eq:dtasd} holds with $\varkappa$ replaced by $\kappa - \lambda$, which once again implies \eqref{eq:dtarea}.

\subsection{The intermediate evolution flow}
In the axisymmetric setting, the intermediate evolution flow in \eqref{eq:imS} can be written as:
\begin{equation}\label{eq:xtSALK}
(\pol x\cdot\pol\ek_1)\,\pol x_t\cdot\pol\nu = 
- \left[\pol x\cdot\pol\ek_1\,\mathcal{Y}_s \right]_s\,,\quad
-\tfrac1\xi\,[\pol x\cdot\pol\ek_1\,\mathcal{Y}_s]_s + \tfrac1\alpha\,\pol
x\cdot\pol\ek_1\,\mathcal{Y} = \pol x\cdot\pol\ek_1\,\varkappa
\qquad \text{in } \mbI \times (0,T]\,,
\end{equation}
together with the boundary conditions \eqref{eq:bcs}
for $t\in(0,T]$, as well as \eqref{eq:sdbdn} with $\varkappa_\rho$ replaced by
$\mathcal{Y}_\rho$ in \eqref{eq:bc5}. 
It is straightforward to adapt the weak formulations 
from \S\ref{sec:sdf} to (\ref{eq:xtSALK}). 
For example, generalizing \eqref{eq:sdstabab} to 
(\ref{eq:xtSALK}) yields the following weak formulation.

\noindent
Let $\pol x(0) \in \Vpartialzero$. For $t \in (0,T]$
find $\pol x(t) \in [H^1(\mbI)]^2$, with $\pol x_t(t) \in \Vpartial$, 
and $(\mathcal{Y}(t), \varkappa(t))\in [H^1(\mbI)]^2$ such that
\begin{subequations} \label{eq:intstababc}
\begin{alignat}{2}
\Big\langle(\pol x\cdot\pol\ek_1)\,\pol x_t\cdot\pol\nu,~\chi\,|\pol x_\rho|\Big\rangle
- \Big\langle\pol x\cdot\pol\ek_1\, \mathcal{Y}_\rho,~ \chi_\rho\,|\pol x_\rho|^{-1} \Big\rangle&=0 &&
\qquad \forall\ \chi \in H^1(\mbI)\,, \label{eq:intstaba} \\
\frac1\xi\, \Big\langle\pol x\cdot\pol\ek_1 \, \mathcal{Y}_\rho,~ \zeta_\rho\,
|\pol x_\rho|^{-1} \Big\rangle
+ \Big\langle\pol x\cdot\pol\ek_1 \left[ \alpha^{-1}\,\mathcal{Y} -
\varkappa \right],~ \zeta\,|\pol x_\rho| \Big\rangle & = 0 &&
\qquad \forall\ \zeta \in H^1(\mbI)\,, \label{eq:intstabb}\\
\Big\langle \pol x\cdot\pol\ek_1\,
\varkappa\,\pol\nu,~\pol\eta\, |\pol x_\rho| \Big\rangle + \Big\langle\pol\ek_1,~\pol\eta \,|\pol x_\rho|\Big\rangle 
+\Big\langle
(\pol x\cdot\pol\ek_1)\,\pol x_\rho,~\pol\eta_\rho|\pol x_\rho|^{-1}\Big\rangle
&= 0 &&\qquad \forall\ \pol\eta \in \Vpartial\,.
\label{eq:intstabc}
\end{alignat}
\end{subequations}
By setting $\chi = 1$ in \eqref{eq:intstaba} and making use of \eqref{eq:dtm}, we obtain \eqref{eq:dtvol}, as before. Choosing $\chi = \mathcal{Y}$ in \eqref{eq:intstaba}, $\zeta = \xi(\varkappa - \alpha^{-1}\mathcal{Y})$ in \eqref{eq:intstabb} and $\pol\eta = \pol x_t$ in \eqref{eq:intstabc}, and using \eqref{eq:dta} yields 
\[
\frac1{2\pi}\,\ddt\,A(\pol x(t)) = -\frac{1}{\alpha}\left\langle\pol x\cdot\pol\ek_1\,|\mathcal{Y}_\rho|^2,|\pol x_\rho|^{-1}\right\rangle -\xi \left\langle\pol x\cdot\pol\ek_1\left|\varkappa - \alpha^{-1}\mathcal{Y}\right|^2, |\pol x_\rho|\right\rangle \leq 0,
\]
which implies that \eqref{eq:dtarea} holds.

For completeness we mention that the weak formulation of (\ref{eq:xtSALK}) 
corresponding to \eqref{eq:sdweak}, i.e.\ with the curve's curvature $\kappa$
being a variable, rather than $\varkappa$, is given 
by (\ref{eq:intstaba}), (\ref{eq:sdweakb}) and (\ref{eq:intstabb}) with
$\varkappa$ replaced by $\kappa - \lambda$. 
 
\subsection{Conserved mean curvature flow} 
The volume preserving mean curvature flow \eqref{eq:cmcfS}, in the axisymmetric
setting, can be formulated as follows.
\begin{equation}\label{eq:xtcmcf}
\pol x_t\cdot\pol\nu = 
\varkappa - \frac{\langle \pol x \cdot \pol\ek_1, \varkappa\,|\pol x_\rho| \rangle}
{\langle  \pol x \cdot \pol\ek_1, |\pol x_\rho| \rangle}
\qquad \text{in } \mbI \times (0,T]\,,
\end{equation}
together with the boundary conditions \eqref{eq:bcs}
for $t\in(0,T]$, as well as \eqref{eq:bc4}, \eqref{eq:bc4a}. 
The weak formulation in the spirit of \eqref{eq:sdstabab} is then given as
follows.

\noindent
Let $\pol x(0) \in \Vpartialzero$. For $t \in (0,T]$
we find $\pol x(t) \in [H^1(\mbI)]^2$, with $\pol x_t(t) \in \Vpartial$, 
and $\varkappa(t)\in L^2(\mbI)$ such that
\begin{subequations} \label{eq:weakcmcf}
\begin{alignat}{2}
 \Bigl\langle (\pol x\cdot\pol\ek_1)\,\pol x_t\cdot\pol\nu,~\chi\,|\pol x_\rho|\Bigr\rangle
- \Big\langle \pol x\cdot\pol\ek_1\,\varkappa 
- \pol x\cdot\pol\ek_1\,\frac{\langle \pol x \cdot \pol\ek_1, \varkappa\,|\pol x_\rho| \rangle}
{\langle  \pol x \cdot \pol\ek_1, |\pol x_\rho| \rangle}
, ~\chi\,|\pol x_\rho| \Big\rangle&=0 &&
\qquad \forall\ \chi \in H^1(\mbI)\,, \label{eq:weakcmcfa} \\
\Big\langle (\pol x\cdot\pol\ek_1)\,
\varkappa\,\pol\nu,~\pol\eta\, |\pol x_\rho| \Big\rangle + \Big\langle\pol\ek_1,~\pol\eta\,|\pol x_\rho| \Big\rangle 
+\Big\langle
(\pol x\cdot\pol\ek_1)\,\pol x_\rho,~\pol\eta_\rho|\pol x_\rho|^{-1}\Big\rangle &=0 && \qquad \forall\ \pol\eta \in \Vpartial\,.
\label{eq:weakcmcfb}
\end{alignat}
\end{subequations}
Choosing $\chi = 1$ in \eqref{eq:weakcmcfa} and recalling 
\eqref{eq:dtm} yields that \eqref{eq:dtvol} is satisfied, 
while choosing $\chi=\varkappa$ in 
\eqref{eq:weakcmcfa} and $\pol\eta = \pol x_t$ in \eqref{eq:weakcmcfb}, 
on noting \eqref{eq:dta}, shows that 
\[
-\frac1{2\pi}\,
\ddt\, A(\pol x(t)) = \langle \pol x \cdot\pol\ek_1,
|\varkappa|^2\,|\pol x_\rho|\rangle
- 
\frac{\left(\langle \pol x\cdot\pol\ek_1,\varkappa\,|\pol x_\rho| \rangle\right)^2}{\langle \pol x\cdot\pol\ek_1,|\pol x_\rho| \rangle}
\geq 0\,,
\]
where we have noted the Cauchy--Schwarz inequality. This proves that also
\eqref{eq:dtarea} holds.

Analogously to before, the weak formulation in the spirit of \eqref{eq:sdweak}
is given by (\ref{eq:weakcmcfa}), with $\varkappa$ replaced by 
$\kappa - \lambda$, and (\ref{eq:sdweakb}). 
 
\setcounter{equation}{0}
\section{Finite element approximations}
\label{sec:fd}

Let $0= t_0 < t_1 < \ldots < t_{M-1} < t_M = T$ be a
partitioning of $[0,T]$ into possibly variable time steps 
$\ttau_m = t_{m+1} - t_{m}$, $m=0,\ldots, M-1$. 
We set $\ttau = \max_{m=0,\ldots, M-1}\ttau_m$. Besides, let $[0,1]=\cup_{j=1}^J \mbI_j$, $J\geq3$, be a
uniform partition of $[0,1]$ into intervals given by the nodes $q_j$,
$\mbI_j=[q_{j-1},q_j]$, 
i.e.\ $q_j = j\,h$ with $h = J^{-1}$, $j=0,\ldots, J$.
We define the finite element spaces as
\[
V^h := \Bigl\{\chi \in C^0(\overline {\mbI}) : \chi\!\mid_{\mbI_j} \nn\
\text{is affine}\ \forall\ j=1,\ldots, J\Bigr\},\qquad
\Vpartialzero^h:= [V^h]^2\cap \Vpartialzero,\qquad \Vpartial^h = [V^h]^2\cap\Vpartial.
\]

Let $\{\pol X^m\}_{0\leq m\leq M}$, with $\pol X^m\in \Vpartialzero^h$, 
be an approximation to $\{\pol x(t)\}_{t\in[0,T]}$ and define
$\Gamma^m = \pol X^m\left(\overline{\mbI}\right)$. Throughout this section we assume
that
\begin{equation*} 
\pol X^m \cdot\pol\ek_1 > 0 \quad \text{in }\
\overline {\mbI}\setminus \partial_0 \mbI
\quad\text{and}\quad
|\pol{X}^m_\rho| > 0 \quad \text{in } \mbI
\qquad 0\leq m\leq M,
\end{equation*}
so that we can set 
\begin{equation*} 
\pol\tau^m = \pol X^m_s = \frac{\pol X^m_\rho}{|\pol X^m_\rho|} 
\qquad \mbox{and} \qquad \pol\nu^m = -(\pol\tau^m)^\perp\,.
\end{equation*}

The main novelty of this paper is the introduction of fully discrete
finite element approximations that satisfy an exact volume preservation for
the three flows \eqref{eq:sdS}, \eqref{eq:imS} and \eqref{eq:cmcfS}. Starting
from the ideas in \cite{Jiang2021,Zhao2021}, it turns out that the crucial
ingredient is the appropriate treatment of the quantity 
$\pol f=\pol x\cdot\pol\ek_1\,|\pol x_\rho|\,\pol\nu$ in the weak formulations
of the flow equations.
Let us define $\pol f^{m+\frac{1}{2}} \in [L^\infty(\mbI)]^2$ via
\begin{align}
\label{eq:averagenor}
\pol f^{m+\frac{1}{2}}= -\frac{1}{6}\Bigl(2(\pol X^m\cdot\pol\ek_1)\,\pol X^m_{\rho}+2(\pol X^{m+1}\cdot\pol\ek_1)\,\pol X^{m+1}_{\rho} + (\pol X^m\cdot\pol\ek_1)\,\pol X^{m+1}_{\rho} + (\pol X^{m+1}\cdot\pol\ek_1)\,\pol X_{\rho}^m\Bigr)^\perp.
\end{align}
Then we can prove the following lemma.
\begin{lemma} \label{lem:MmM}
Let $\pol X^m \in \Vhpartialzero$ and $\pol X^{m+1} \in \Vhpartialzero$. Then
it holds that
\begin{equation} \label{eq:MMf}
M(\pol X^{m+1})-M(\pol X^m) = 2\pi\,\Bigl\langle\pol X^{m+1}-\pol X^m,~\pol f^{m+\frac{1}{2}}\Bigr\rangle.
\end{equation}
\end{lemma}
\begin{proof}
For $\alpha\in[0,1]$,
we define the finite element functions $\pol X^h(\alpha)\in\Vhpartialzero$ via linear interpolation between $\pol X^m$ and $\pol X^{m+1}$:
\begin{align}
\label{eq:xlinear}
\pol X^h(\rho,\alpha):=(1-\alpha)\pol X^m(\rho) + \alpha\pol X^{m+1}(\rho),\quad \rho\in \mbI, \quad 0\leq \alpha\leq 1.
\end{align}
Therefore we have $\pol X^h(\cdot, 0)=\pol X^m$, $\pol X^h(\cdot,1)=\pol X^{m+1}$ and $\pol X^h_\alpha(\cdot, \alpha) = \pol X^{m+1}-\pol X^m$. 
For the family of surfaces generated by $\{\pol X(\alpha)\}_{\alpha \in [0,1]}$
it holds, similarly to \eqref{eq:dtm}, that
\begin{align}
\dd\alpha\,M(\pol X^h(\alpha)) =2\pi \left\langle \pol X^h(\alpha)\cdot\pol\ek_1, \pol X^h_{\alpha}(\alpha)\cdot\pol\nu^h(\alpha)\,|\pol X^h_{\rho}(\alpha)| \right\rangle,
\label{eq:masspro}
\end{align}
where $\pol \nu^h(\cdot,\alpha)=-|\pol X^h_\rho(\cdot,\alpha)|^{-1}(\pol X^h_\rho(\cdot,\alpha))^\perp$,
compare also with Theorem 71 in \cite{Barrett20}. 
Integrating \eqref{eq:masspro} from 0 to 1 with respect to $\alpha$ yields
\begin{align}
M(\pol X^{m+1}) - M(\pol X^m) &=2\pi\int_0^1\int_{\mbI}(\pol X^h(\alpha)\cdot\pol e_1)\left(\pol X^{m+1}-\pol X^m\right)\cdot(-\pol X^h_\rho(\alpha))^\perp\drho\rd\alpha \nn\\
&=2\pi\int_{\mbI}(\pol X^{m+1}-\pol X^m)\cdot\int_0^1(\pol X^h(\alpha)\cdot\pol\ek_1)(-\pol X^h_{\rho}(\alpha))^\perp\rd\alpha\drho\nn.
\end{align}
By \eqref{eq:xlinear}, we note that $(\pol X^h(\alpha)\cdot\pol\ek_1)(-\pol X^h_{\rho}(\alpha))^\perp$ is a quadratic function of $\alpha$. Applying Simpson's rule and using \eqref{eq:averagenor} yields
\begin{align}
\int_0^1(\pol X^h\cdot\pol\ek_1)(-\pol X^h_{\rho})^\perp\rd\alpha &= -\frac{1}{6}\left((\pol X^m\cdot\pol\ek_1)\pol X_\rho^m + 4(\pol X^{m+\frac{1}{2}}\cdot\pol\ek_1)\pol X^{m+\frac{1}{2}}_\rho + (\pol X^{m+1}\cdot\pol\ek_1)\pol X^{m+1}_\rho\right)^\perp\nn\\
&=-\frac{1}{6}\left(2(\pol X^m\cdot\pol\ek_1)\,\pol X^m_{\rho}+2(\pol X^{m+1}\cdot\pol\ek_1)\,\pol X^{m+1}_{\rho} + (\pol X^m\cdot\pol\ek_1)\,\pol X^{m+1}_{\rho} + (\pol X^{m+1}\cdot\pol\ek_1)\,\pol X_{\rho}^m\right)^\perp=\pol f^{m+\frac{1}{2}},
\end{align}
where we have denoted $\pol X^{m+\frac{1}{2}} = \frac{1}{2}(\pol X^m + \pol X^{m+1})$. Therefore we have \eqref{eq:MMf} as claimed.
\end{proof}

\subsection{For the surface diffusion flow}
For the weak formulation \eqref{eq:sdstabab}, we propose the following discretized scheme. Let $\pol X^0 \in \Vhpartialzero$. For $m=0,\ldots,M-1$, 
find $(\delta\pol X^{m+1}, \varkappa^{m+1}) \in \Vhpartial \times V^h$,
where $\pol X^{m+1} = \pol X^m + \delta \pol X^{m+1}$, 
such that
\begin{subequations}  \label{eq:sdstabfd}
\begin{alignat}{2}
\frac{1}{\ttau_m}\left\langle\pol X^{m+1} - \pol X^m,
~\chi\,\pol f^{m+\frac{1}{2}}\right\rangle
- \left\langle\pol X^m\cdot\pol\ek_1 \,
\varkappa^{m+1} _\rho, 
~\chi_\rho\,|\pol X^m_\rho|^{-1}\right\rangle&=0&&
\qquad \forall\ \chi \in V^h\,, \label{eq:sdstabfda}\\
\left\langle
\varkappa^{m+1}\,\pol f^{m+\frac{1}{2}}, ~\pol\eta\right\rangle
+ \left\langle \pol\eta \cdot\pol\ek_1, ~|\pol X^{m+1}_\rho|\right\rangle
+ \left\langle (\pol X^m\cdot\pol\ek_1)\,
\pol X^{m+1}_\rho,~\pol\eta_\rho\, |\pol X^m_\rho|^{-1} \right\rangle
&= 0 &&\qquad \forall\ \pol\eta \in \Vhpartial\,.
\label{eq:sdstabfdb}
\end{alignat}
\end{subequations}
We note that the scheme \eqref{eq:sdstabfd} is very close to \cite[(4.3)]{Barrett19}. The main difference is that there the combined term $\pol f=\pol x\cdot\pol\ek_1\,|\pol x_\rho|\,\pol\nu$ is approximated fully explicitly. Here in \eqref{eq:sdstabfda} we choose the semi-implicit treatment $\pol f^{m+\frac{1}{2}}$ from \eqref{eq:averagenor} to obtain exact volume preservation. In addition, we also employ $\pol f^{m+\frac{1}{2}}$ in \eqref{eq:sdstabfdb} to obtain an unconditionally stable method. 

We have the following theorem for the discretization \eqref{eq:sdstabfd}, which mimics the energy dissipation and the volume preservation on the discrete level.
\begin{thm}[stability and volume conservation]\label{thm:sfvc}
Let $(\pol X^{m+1},\kappa^{m+1})$ be a solution to \eqref{eq:sdstabfd}.
Then it holds that
\begin{align}
\label{eq:stabe}
A(\pol X^{m+1}) + 2\pi\,\ttau_m\Bigl\langle\pol X^m\cdot\pol\ek_1\,\varkappa_{\rho}^{m+1},~\varkappa_{\rho}^{m+1}|\pol X^m_{\rho}|^{-1}\Bigr\rangle \leq A(\pol X^m).
\end{align}
Moreover, it holds that
\begin{equation} \label{eq:stabm}
M(\pol X^{m+1})=M(\pol X^m).
\end{equation}
\end{thm}
\begin{proof}
Choosing $\chi = \varkappa^{m+1}$ in \eqref{eq:sdstabfda} and $\pol\eta = \pol X^{m+1} - \pol X^m$ in \eqref{eq:sdstabfdb}, and then combining the two equations, we obtain
\begin{align}\label{eq:eime}
&\ttau_m\Bigl\langle\pol X^m\cdot\pol\ek_1\,\varkappa_{\rho}^{m+1},~\varkappa_{\rho}^{m+1}|\pol X^m_{\rho}|^{-1}\Bigr\rangle+\left\langle (\pol X^{m+1}-\pol X^m) \cdot\pol\ek_1, ~|\pol X^{m+1}_\rho|\right\rangle
+ \left\langle (\pol X^m\cdot\pol\ek_1)\,
\pol X^{m+1}_\rho,~(\pol X^{m+1} - \pol X^m)_\rho\, |\pol X^m_\rho|^{-1} \right\rangle =0.
\end{align}
Using the inequality $\pol a\cdot(\pol a - \pol b)\geq |\pol b|(|\pol a| - |\pol b|)$ for $\pol a, \pol b \in \bR^2$, we have
\begin{align}
&\left\langle (\pol X^{m+1}-\pol X^m) \cdot\pol\ek_1, ~|\pol X^{m+1}_\rho|\right\rangle
+ \left\langle (\pol X^m\cdot\pol\ek_1)\,
\pol X^{m+1}_\rho,~(\pol X^{m+1} - \pol X^m)_\rho\, |\pol X^m_\rho|^{-1} \right\rangle\nn\\
&\hspace{2cm}\geq \left\langle(\pol X^{m+1}-\pol X^m) \cdot\pol\ek_1, ~|\pol X^{m+1}_\rho|\right\rangle + \left\langle\pol X^m\cdot\pol\ek_1 ,~|\pol X^{m+1}_\rho| - |\pol X^m_{\rho}| \right\rangle \nn\\
&\hspace{2cm}= \Bigl\langle\pol X^{m+1}\cdot\pol\ek_1,~|\pol X^{m+1}_{\rho}|\Bigr\rangle - \Bigl\langle\pol X^m\cdot\pol\ek_1,~|\pol X^m_\rho|\Bigr\rangle= \frac{1}{2\pi}\left(A(\pol X^{m+1})-A(\pol X^m)\right),
\label{eq:sfstableimed1}
\end{align}
on recalling \eqref{eq:TE}. Combining \eqref{eq:eime} and \eqref{eq:sfstableimed1} 
yields the desired stability result \eqref{eq:stabe}. 

Moreover, setting $\chi = 1$ in \eqref{eq:sdstabfda} yields
\[
\Bigl\langle\pol X^{m+1}-\pol X^m,~\pol f^{m+\frac{1}{2}}\Bigr\rangle=0,
\]
which implies \eqref{eq:stabm} thanks to Lemma~\ref{lem:MmM}.
\end{proof}

We next consider the discretization of the weak formulation in 
\eqref{eq:sdweak}, leading to a scheme with an equidistribution property. 
Following the approach from \cite[(3.8)]{Barrett19}, we first note from (\ref{eq:bc4}) that
\begin{equation*} 
\lim_{\rho\to \rho_{_0}}\lambda(\rho,~t)=\lim_{\rho\to \rho_{_0}}
\frac{\pol\nu(\rho,t)\cdot\pol\ek_1}{\pol x(\rho,t)\cdot\pol\ek_1} 
= \lim_{\rho\to \rho_{_0}} 
\frac{\pol\nu_\rho(\rho,t)\cdot\pol\ek_1}{\pol x_\rho(\rho,t)\cdot\pol\ek_1}
= \pol\nu_s(\rho_{_0},t)\cdot\pol\tau(\rho_{_0},t) 
= -\kappa(\rho_{_0},t)
\quad \forall\ \rho_{_0}\in\partial_0\mbI\,,\
\forall\ t \in [0,T]\,. 
\end{equation*}
Thus, for $\kappa^{m+1} \in V^h$, we introduce $\lambda^{m+\frac{1}{2}}(\kappa^{m+1}) \in V^h$  to avoid the degeneracy in the discretization such that
\begin{equation*} 
[\lambda^{m+\frac{1}{2}}(\kappa^{m+1})](q_j) = \begin{cases}
- \kappa^{m+1}(q_j) & q_j \in \partial_0\mbI\,,\\[0.7em]
\dfrac{\pol\omega^m(q_j)\cdot\pol\ek_1}{\pol X^m(q_j)\cdot\pol\ek_1}
& {\text{otherwise}}\,, 
\end{cases}
\end{equation*}
where $\pol\omega^m \in [V^h]^2$ is the mass-lumped 
$L^2$--projection of $\pol\nu^m$ onto $[V^h]^2$, i.e.\
\begin{equation*} 
\left\langle\pol\omega^m, \pol\varphi \, |\pol X^m_\rho| \right\rangle^h 
= \left\langle \pol\nu^m, \pol\varphi \, |\pol X^m_\rho| \right\rangle
\qquad \forall\ \pol\varphi\in[V^h]^2\,.
\end{equation*}
Here $\langle\cdot,\cdot\rangle^h$ denotes the usual 
mass lumped $L^2$--inner product, which 
for two piecewise continuous functions $\pol v, \pol w$, 
with possible jumps at the nodes $\{q_j\}_{j=1}^J$, is defined via
\[
\Big\langle \pol v, \pol w \Big\rangle^h = \tfrac12\,h\sum_{j=1}^J 
\left[(\pol v\cdot\pol w)(q_j^-) + (\pol v \cdot \pol w)(q_{j-1}^+)\right],
\]
where
$g(q_j^\pm)=\underset{\delta\searrow 0}{\lim}\ g(q_j\pm\delta)$.

Then for the weak formulation \eqref{eq:sdweak} we propose the following finite element approximation. Let $\pol X^0 \in \Vhpartialzero$. For $m=0,\ldots,M-1$, 
find $(\delta\pol X^{m+1}, \kappa^{m+1}) \in \Vhpartial \times V^h$,
where $\pol X^{m+1} = \pol X^m + \delta \pol X^{m+1}$, 
such that
\begin{subequations} \label{eq:sdfd}
\begin{alignat}{2}
\dfrac{1}{\ttau_m}\left\langle\pol X^{m+1} - \pol X^m,~\chi\,\pol f^{m+\frac{1}{2}}\right\rangle- \left\langle \pol X^m\cdot\pol\ek_1 \left[
\kappa^{m+1} - \lambda^{m+\frac{1}{2}}(\kappa^{m+1})\right]_\rho,~
\chi_\rho\,|\pol X^m_\rho|^{-1}\right\rangle&=0&&
\qquad \forall\ \chi \in V^h\,, \label{eq:sdfda}
\\
\Big\langle\kappa^{m+1}\,\pol\nu^m, ~\pol\eta\,|\pol X^m_\rho|\Big\rangle^h
+ \Big\langle\pol X^{m+1}_\rho, ~\pol\eta_\rho\,|\pol X^m_\rho|^{-1}\Big\rangle 
&= 0 &&\qquad \forall\ \pol\eta \in \Vhpartial\,.
\label{eq:sdfdb}
\end{alignat}
\end{subequations}
It does not appear possible to recover the energy stability for the discretized scheme \eqref{eq:sdfd}. Nevertheless, in a similar manner to \eqref{eq:sdstaba}, the approximation using $\pol f^{m+\frac{1}{2}}$ in the first term of Eq.~\eqref{eq:sdfda} contributes to the property of volume conservation. Besides, the first term in \eqref{eq:sdfdb} is approximated by using the mass lumped inner product, which leads to the property of equidistribution, i.e., the mesh points on $\Gamma^m$ tend to be distributed at evenly spaced arc length. For a semidiscrete approximation this can be made rigorous, see Remark~3.1 in \cite{Barrett19}. A detailed discussion of this equidistribution property can be found in \cite{Barrett07,Barrett20}. 
 
It is straightforward to show that the discretized scheme \eqref{eq:sdfd} satisfies the exact volume conservation. 
\begin{thm}[volume conservation] \label{thm:sdv}
Let $(\pol X^{m+1},\kappa^{m+1})$ be a solution to \eqref{eq:sdfd}.
Then it holds that
\begin{equation*} 
M(\pol X^{m+1})=M(\pol X^m).
\end{equation*}
\end{thm}
\begin{proof}
Setting $\chi=1$ in \eqref{eq:sdfda} and noting Lemma~\ref{lem:MmM} yields the
desired result. 
\end{proof}

\subsection{For the intermediate evolution flow}
For the weak formulation \eqref{eq:intstababc}, we propose the following fully discretized scheme: Let $\pol X^0 \in \Vhpartialzero$. For $m=0,\ldots,M-1$, 
find $(\delta\pol X^{m+1}, \mathcal{Y}^{m+1}, \varkappa^{m+1}) 
\in \Vhpartial \times [V^h]^2$, 
where $\pol X^{m+1} = \pol X^m + \delta \pol X^{m+1}$, such that
\begin{subequations}  \label{eq:stabimabc}
\begin{alignat}{2}
\frac{1}{\ttau_m}\Big\langle\pol X^{m+1} - \pol X^m,
~\chi\,\pol f^{m+\frac{1}{2}}\,\Big\rangle -\Big\langle\pol X^m\cdot\pol\ek_1 \, \mathcal{Y}^{m+1}_\rho, 
~\chi_\rho\,|\pol X^m_\rho|^{-1}\Big\rangle&=0&&
\qquad \forall\ \chi \in V^h\,, \label{eq:stabima}\\
\tfrac1\xi \Big\langle\pol X^m\cdot\pol\ek_1 \, \mathcal{Y}^{m+1}_\rho, 
~\zeta_\rho\,|\pol X^m_\rho|^{-1}\Big\rangle
+ \Big\langle\pol X^m\cdot\pol\ek_1 \left[ \alpha^{-1}\,\mathcal{Y}^{m+1} -
\varkappa^{m+1} \right], 
~\zeta\,|\pol X^m_\rho|\Big\rangle & = 0&&
\qquad \forall\ \zeta \in V^h\,, \label{eq:stabimb}\\
\left\langle
\varkappa^{m+1}\,\pol f^{m+\frac{1}{2}}, ~\pol\eta\right\rangle
+ \left\langle \pol\eta \cdot\pol\ek_1, ~|\pol X^{m+1}_\rho|\right\rangle
+ \left\langle (\pol X^m\cdot\pol\ek_1)\,
\pol X^{m+1}_\rho,~\pol\eta_\rho\, |\pol X^m_\rho|^{-1} \right\rangle
& = 0 &&\qquad \forall\ \pol\eta \in \Vhpartial\,.
\label{eq:stabimc}
\end{alignat}
\end{subequations}
For the discretization in \eqref{eq:stabimabc} we can prove unconditional energy decay and exact volume conservation.
\begin{thm}[stability and volume conservation] \label{thm:imstabv}
Let $(\pol X^{m+1},\mathcal{Y}^{m+1},\kappa^{m+1})$ be a solution to \eqref{eq:stabimabc}. Then it holds that
\begin{align}
A(\pol X^{m+1}) + 
\frac{2\pi\,\ttau_m}{\alpha}
\Big\langle\pol X^m\cdot\pol\ek_1 
\,\left|\mathcal{Y}^{m+1}_\rho\right|^2, ~|\pol X^m_\rho|^{-1}\Big\rangle
+ 2\pi\,\ttau_m\,\xi \Big\langle \pol X^m\cdot\pol\ek_1\,
\left|\varkappa^{m+1} - \tfrac1\alpha\,\mathcal{Y}^{m+1}\right|^2 , ~|\pol X^m_\rho| 
\Big\rangle \leq A(\pol X^m)\,.
\label{eq:stabimE}
\end{align}
Moreover, it holds that
\begin{equation} \label{eq:stabimfv}
M(\pol X^{m+1})=M(\pol X^m).
\end{equation}
\end{thm}
\begin{proof}
Setting $\chi = \ttau_m\,\varkappa^{m+1}$ in
(\ref{eq:stabima}),
$\zeta=\ttau_m\,\xi(\varkappa^{m+1}-\frac{1}{\alpha}\mathcal{Y}^{m+1})$ in (\ref{eq:stabimb}) 
and $\pol\eta = \pol X^{m+1} - \pol X^m$ 
in (\ref{eq:stabimc}), and then combining these three equations yields
\begin{align}
&\frac{\Delta t_m}{\alpha}\Big\langle\pol X^m\cdot\pol\ek_1 
\,\left|\mathcal{Y}^{m+1}_\rho\right|^2, |\pol X^m_\rho|^{-1}\Big\rangle + \Delta t_m\,\xi \Big\langle \pol X^m\cdot\pol\ek_1\,
\left|\varkappa^{m+1} - \tfrac1\alpha\,\mathcal{Y}^{m+1}\right|^2 , |\pol X^m_\rho|\Big\rangle +\Big\langle(\pol X^{m+1} - \pol X^m)\cdot\pol\ek_1,~|\pol X^{m+1}_\rho|\Big\rangle\nn\\&\hspace{1cm} + \Big\langle (\pol X^m\cdot\pol\ek_1)\,
\pol X^{m+1}_\rho,~(\pol X^{m+1} - \pol X^m)_\rho\, |\pol X^m_\rho|^{-1} \Big\rangle=0.
\label{eq:imstableimed}
\end{align}
Using \eqref{eq:sfstableimed1} we can recast \eqref{eq:imstableimed} as
\begin{align}
\frac{\Delta t_m}{\alpha}\left(\pol X^m\cdot\pol\ek_1 
\,\left|\mathcal{Y}^{m+1}_\rho\right|^2, |\pol X^m_\rho|^{-1}\right) + \Delta t_m\,\xi \left( \pol X^m\cdot\pol\ek_1\,
\left|\varkappa^{m+1} - \tfrac1\alpha\,\mathcal{Y}^{m+1}\right|^2 , |\pol X^m_\rho|\right) \leq -\frac{1}{2\pi}\left(A(\pol X^{m+1})- A(\pol X^m)\right),\nn
\end{align}
which implies \eqref{eq:stabimE} immediately.

Finally, as before, setting $\chi = 1$ in \eqref{eq:stabimc} and noting 
Lemma~\ref{lem:MmM} yields \eqref{eq:stabimfv}.
\end{proof}

It is not difficult to also introduce the natural finite element approximation 
of the intermediate evolution law \eqref{eq:xtSALK} in the spirit of the 
scheme \eqref{eq:sdfd}.
In fact, that method would be given by \eqref{eq:stabima}, \eqref{eq:stabimb} 
with $\varkappa^{m+1}$, replaced by 
$\kappa^{m+1} - \lambda^{m+\frac12}(\kappa^{m+1})$,
and \eqref{eq:sdfdb}. That new scheme would again satisfy the equidistribution
property induced by \eqref{eq:sdfdb}, together with the exact volume
conservation as stated in Theorem~\ref{thm:sdv}.

\subsection{For the conserved mean curvature flow} 
For the weak formulation \eqref{eq:weakcmcf}, we propose the following fully discrete approximation. Let $\pol X^0 \in \Vhpartialzero$. For $m=0,\ldots,M-1$, 
find $(\delta\pol X^{m+1}, \varkappa^{m+1}) \in \Vhpartial \times V^h$, 
where $\pol X^{m+1} = \pol X^m + \delta \pol X^{m+1}$, such that
\begin{subequations}  \label{eq:stabcmcf}
\begin{align}
&
\frac{1}{\ttau_m}\Big\langle\pol X^{m+1} - \pol X^m,
~\chi\,\pol f^{m+\frac{1}{2}}\,\Big\rangle -\Big\langle\pol X^m\cdot\pol\ek_1 \, \varkappa^{m+1}, ~\chi\,|\pol X^m_\rho|\Big\rangle
+ \frac{\left\langle\pol X^m\cdot\pol\ek_1, \varkappa^{m+1}\,
|\pol X^m_\rho|\right\rangle}
{\left\langle\pol X^m\cdot\pol\ek_1, |\pol X^m_\rho|\right\rangle}
\left\langle \pol X^m\cdot\pol\ek_1, \chi\,|\pol X^m_\rho|\right\rangle
=0
\qquad \forall\ \chi \in V^h\,, \label{eq:stabcmcfa}\\
&
\left\langle
\varkappa^{m+1}\,\pol f^{m+\frac{1}{2}}, ~\pol\eta\right\rangle
+ \left\langle \pol\eta \cdot\pol\ek_1, ~|\pol X^{m+1}_\rho|\right\rangle
+ \left\langle (\pol X^m\cdot\pol\ek_1)\,
\pol X^{m+1}_\rho,~\pol\eta_\rho\, |\pol X^m_\rho|^{-1} \right\rangle
= 0 \qquad \forall\ \pol\eta \in \Vhpartial\,.
\label{eq:stabcmcfb}
\end{align}
\end{subequations}
Similarly to before, we can prove unconditional energy decay and exact volume conservation for the scheme \eqref{eq:stabcmcf}.
\begin{thm}[stability and volume conservation] 
Let $(\pol X^{m+1},\kappa^{m+1})$ be a solution to \eqref{eq:stabcmcf}. Then it holds that
\begin{equation} \label{eq:stabcmcfA}
\frac1{2\pi}\,A(\pol X^{m}) - \frac1{2\pi}\,A(\pol X^{m+1})
\geq \ttau_m \left\langle\pol X^m\cdot\pol\ek_1\,|\varkappa^{m+1}|^2,
|\pol X^m_\rho|\right\rangle
- \ttau_m
\frac{\left(\left\langle \pol X^m\cdot\pol\ek_1,\varkappa^{m+1}\,
|\pol X^m_\rho| \right\rangle\right)^2}
{\left\langle\pol X^m\cdot\pol\ek_1,|\pol X^m_\rho| \right\rangle} \geq 0\,.
\end{equation}
Moreover, it holds that
\begin{equation} \label{eq:stabcmcfv}
M(\pol X^{m+1})=M(\pol X^m).
\end{equation}
\end{thm}
\begin{proof}
Setting $\chi = \ttau_m\,\varkappa^{m+1}$ in (\ref{eq:stabcmcfa})
and $\pol\eta = \pol X^{m+1} - \pol X^m$ in (\ref{eq:stabcmcfb}),
combining and recalling \eqref{eq:sfstableimed1} yields 
\eqref{eq:stabcmcfA}, where we also note a Cauchy--Schwarz inequality.
In addition, as before, choosing $\chi = 1$ in \eqref{eq:stabcmcfa} 
and noting Lemma~\ref{lem:MmM} yields \eqref{eq:stabcmcfv}.
\end{proof}

We can also introduce the natural finite element approximation 
of \eqref{eq:xtcmcf} in the spirit of the scheme \eqref{eq:sdfd}.
The new method would be given by \eqref{eq:stabcmcfa},
with $\varkappa^{m+1}$ replaced by 
$\kappa^{m+1} - \lambda^{m+\frac12}(\kappa^{m+1})$,
and \eqref{eq:sdfdb}. That new scheme would again satisfy the equidistribution
property induced by \eqref{eq:sdfdb}, together with the exact volume
conservation as stated in Theorem~\ref{thm:sdv}.

\subsection{The nonlinear solver}
We note that all our introduced schemes lead to systems of nonlinear equations, which can be solved with the help of a Newton's method. For example, for \eqref{eq:sdstabfd},  given the initial guess $\left(\pol X^{m+1,0}, \varkappa^{m+1,0}\right)\in\Vpartialzero \times V^h$, for $i\geq 0$ we seek the Newton's direction $\Bigl(\pol X^{\delta},~\varkappa^{\delta}\Bigr)\in\Vpartial^h \times V^h$ such that the following equations hold
\begin{subequations}
\label{eq:sdstabfdnewton}
\begin{align}
&\frac{1}{\ttau_m}\Big\langle\pol X^{\delta},~\chi\,\pol f^{m+\frac{1}{2},i}\Big\rangle + \frac{1}{\ttau_m}\Bigl\langle\pol X^{m+1, i} - \pol X^m,
~\chi\,\pol f^{m+\frac{1}{2}, i}_{\delta}\Big\rangle-\Big\langle\pol X^m\cdot\pol\ek_1 \,
\varkappa^{\delta} _\rho, 
~\chi_\rho\,|\pol X^m_\rho|^{-1}\Big\rangle\nn \\&\hspace{1cm}=-\frac{1}{\ttau_m}\Big\langle\pol X^{m+1, i} - \pol X^m,
~\chi\,\pol f^{m+\frac{1}{2}, i}\Big\rangle
+ \Big\langle\pol X^m\cdot\pol\ek_1 \,
\varkappa^{m+1, i} _\rho, 
~\chi_\rho\,|\pol X^m_\rho|^{-1}\Big\rangle\qquad \forall\ \chi \in V^h,\\
&\Big\langle
\varkappa^{\delta}\,\pol f^{m+\frac{1}{2}, i}, ~\pol\eta\Big\rangle + \Big\langle
\varkappa^{m+1, i}\,\pol f^{m+\frac{1}{2}}_{\delta}, ~\pol\eta\Big\rangle+\Big\langle \pol\eta \cdot\pol\ek_1, ~|\pol X^{m+1,i}_\rho|^{-1}\pol X^{m+1,i}_\rho\cdot\pol X^\delta_\rho\Big\rangle+\Big\langle (\pol X^m\cdot\pol\ek_1)\,
\pol X^{\delta}_\rho,~\pol\eta_\rho\, |\pol X^m_\rho|^{-1} \Big\rangle\nn\\
&\hspace{1cm} =-\Bigl\langle
\varkappa^{m+1, i}\,\pol f^{m+\frac{1}{2},i}, ~\pol\eta\Big\rangle
- \Big\langle\pol\eta \cdot\pol\ek_1, ~|\pol X^{m+1,i}_\rho|\Big\rangle
- \Big\langle (\pol X^m\cdot\pol\ek_1)\,
\pol X^{m+1, i}_\rho,~\pol\eta_\rho\, |\pol X^m_\rho|^{-1} \Big\rangle
\qquad \forall\ \pol\eta \in \Vpartial^h,
\end{align}
\end{subequations}
where $\pol f^{m+\frac{1}{2},i}\in [L^\infty(\mbI)]^2$ and $\pol f^{m+\frac{1}{2},i}_\delta\in [L^\infty(\mbI)]^2$ are given by
\begin{align*}
\pol f^{m+\frac{1}{2},i} &:=-\frac{1}{6}\Bigl[2(\pol X^m\cdot\pol\ek_1)\,\pol X^m_{\rho}+2(\pol X^{m+1,i}\cdot\pol\ek_1)\,\pol X^{m+1,i}_{\rho} + (\pol X^m\cdot\pol\ek_1)\,\pol X^{m+1,i}_{\rho} + (\pol X^{m+1,i}\cdot\pol\ek_1)\,\pol X_{\rho}^m\Bigr]^\perp,\\
\pol f^{m+\frac{1}{2},i}_\delta:&=-\frac{1}{6}\left[(\pol X^\delta\cdot\pol\ek_1)\,(2 \pol X^{m+1,i}_\rho+\pol X_\rho^m) + (2\pol X^{m+1,i}\cdot\pol\ek_1 + \pol X^m\cdot\pol\ek_1)\,\pol X^\delta_\rho\right]^\perp.
\end{align*}
We then update 
\begin{align} \label{eq:newtonup}
\pol X^{m+1,i+1} = \pol X^{m+1,i} + \pol X^\delta,\qquad \varkappa^{m+1,i+1} = \varkappa^{m+1,i}+\varkappa^\delta.
\end{align}
For each $m\geq0$, we can choose the initial guess $\pol X^{m+1,0}=\pol X^m$, $\varkappa^{m+1,0}=\varkappa^m$, and then repeat the iterations in \eqref{eq:sdstabfdnewton}, \eqref{eq:newtonup} until the following conditions hold:

\begin{align}
&\vnorm{\pol X^{m+1,i+1}-\pol X^{m+1,i}}_{\infty}=\max_{1\leq j\leq J}|\pol X^{m+1,i+1}(q_j) - \pol X^{m+1,i}(q_j)|\leq{\rm tol},\nn\\
&\vnorm{\varkappa^{m+1,i+1}-\varkappa^{m+1,i}}_{\infty}=\max_{1\leq j\leq J}|\varkappa^{m+1,i+1}(q_j) - \varkappa^{m+1,i}(q_j)|\leq{\rm tol},\label{eq:stop}
\end{align}
where ${\rm tol}$ is the chosen tolerance. 

\subsection{Extension to non-neutral external boundaries}
So far, to simplify the presentation, we have only considered the case of
neutral external boundaries, when the axisymmetric surfaces $\mathcal{S}(t)$
have a nonempty boundary. This then leads to the $90^\circ$ contact angle
conditions induced by \eqref{eq:bc4a}. In physical applications, however, it
can be of interest to also consider scenarios, where a difference in contact
energy densities between the external boundary and the interior and exterior
of $\mathcal{S}(t)$, respectively, leads to other contact angles.
To this end, we consider contact energy contributions as discussed in
\cite{Finn86}, see also \cite[(2.21)]{Barrett10}. 
In the axisymmetric setting, the relevant energy is then given by
\begin{align}
 E(\pol x(t))& = A(\pol x(t)) + A_{\partial_1}(\pol x(t)) 
+ A_{\partial_2}(\pol x(t))\nn\\
 &=A(\pol x(t)) + 2\pi\,\sum_{p\in\partial_1 \mbI} 
\sliprho^{(p)}\,(\pol x(p,t)\cdot\pol\ek_1)\,\pol x(p,t)\cdot\pol\ek_2
+ \pi\,\sum_{p\in\partial_2\mbI} 
\sliprho^{(p)}\,(\pol x(p,t)\cdot\pol\ek_1)^2\,,
 \label{eq:E}
\end{align}
where we recall that $A(\pol x)$ is the surface area of $\mathcal{S}(t)$, 
and $A_{\partial_i}(\pol x)$ represent contact energies of the boundary
$\partial\mathcal{S}(t)$ with an external boundary.
In particular, $\sliprho^{(p)}$, for $p \in \partial_1 \mbI$, 
denotes the change in contact energy
density in the direction of $-\pol\ek_2$, that the two phases separated by
the interface $\mathcal{S}(t)$ have with the infinite cylinder at
the boundary circle of $\mathcal{S}(t)$ represented by $\pol x(p,t)$,
compare with the simulation in Fig.~\ref{fig:icylinder} below.
Similarly, $\sliprho^{(p)}$, for $p \in \partial_2\mbI$,
denotes the change in contact energy
density in the direction of $-\pol\ek_1$, that the two phases separated by
the interface $\mathcal{S}(t)$ have with the hyperplane 
$\bR\times\{\pol x(p,0) \cdot \pol\ek_2\}\times\bR$ 
at the boundary circle of $\mathcal{S}(t)$ represented by
$\pol x(p,t)$,
compare with the simulation in Fig.~\ref{fig:substrate} below.
The contact energy contributions in \eqref{eq:E} lead to the natural contact angle conditions
\begin{subequations} \label{eq:mcbc}
\begin{align} 
(-1)^p\,\pol\tau(p,t)\cdot\pol\ek_2 &= \sliprho^{(p)} 
\qquad p \in \partial_1\mbI\,,\label{eq:mcbc1} \\
(-1)^p\,\pol\tau(p,t)\cdot\pol\ek_1 &= \sliprho^{(p)} 
\qquad p \in \partial_2\mbI\,,\label{eq:mcbc2}
\end{align}
\end{subequations}
for all $t \in (0,T]$. In the case of a neutral external boundary we have $\sliprho^{(0)}=\sliprho^{(1)}=0$, which leads to the $90^\circ$ contact angle conditions in (\ref{eq:bc4a}), and means that (\ref{eq:E}) reduces to $E(\pol x) = A(\pol x)$. We refer to \cite{Finn86, Barrett10, Jiang19a} for more details on contact angles and contact energies. We note that in \cite{Bao17,Jiang19a,Zhao20} there are discussions of the relaxed contact angle conditions, which give rise to dynamic contact angles that may be different from the equilibrium contact angles. The treatment of these relaxed conditions is not much different to that of \eqref{eq:mcbc}. Thus we will simply focus on the natural contact angle conditions \eqref{eq:mcbc} in this paper.

Taking the time derivative of $E(\pol x)$, and noting the boundary conditions \eqref{eq:freeslipbc}, yields
\begin{align}
\ddt\, E(\pol x(t)) & = 2\pi\,\int_\mbI \left[\pol x_t\cdot\pol\ek_1
+ \pol x\cdot\pol\ek_1\,\frac{(\pol x_t)_\rho\cdot\pol x_\rho}{|\pol x_\rho|^2}
\right] |\pol x_\rho| \drho 
+ 2\pi\sum_{p\in\partial_1\mbI} 
\sliprho^{(p)}(\pol x(p,t)\cdot\pol\ek_1)\,\pol x_t(p,t)\cdot\pol\ek_2
\nonumber \\ & \qquad
+ 2\pi\sum_{p\in\partial_2\mbI} 
\sliprho^{(p)}\,(\pol x(p,t)\cdot\pol\ek_1)\,\pol x_t(p,t)\cdot\pol\ek_1
\label{eq:dEdt}
\end{align}
as the generalization of \eqref{eq:dta}. Adjusting the weak formulations in 
Section~\ref{sec:wf} is now straightforward. In particular, on the right hand
sides of \eqref{eq:sdstabb}, \eqref{eq:intstabc} and \eqref{eq:weakcmcfb} 
we add the terms
\begin{equation} \label{eq:rhoterms}
 - \sum_{i=1}^2 \sum_{p \in \partial_i\mbI} \sliprho^{(p)}\,
(\pol x(p,t)\cdot\pol\ek_1)\,\pol\eta(p)\cdot\pol\ek_{3-i},
\end{equation}
while to the right hand side of \eqref{eq:sdweakb} we add the terms
\begin{equation} \label{eq:rhoterms2}
- \sum_{i=1}^2 
\sum_{p \in \partial_i\mbI} \sliprho^{(p)}\,\pol\eta(p)\cdot\pol\ek_{3-i}\,.
\end{equation}
In each case, the boundary conditions \eqref{eq:mcbc} will then be weakly
enforced.
Moreover, these new weak formulations still satisfy the volume conservation and energy dissipation properties. For example, using the same testing procedure as before, and noting \eqref{eq:dEdt}, it is straightforward to show that solutions to the adapted \eqref{eq:sdstabab} satisfy \eqref{eq:dtasd} with $A(\pol x(t))$ replaced by $E(\pol x(t))$.

Similarly, the adaptation of our numerical schemes to the more general contact
angles is also straightforward. For the equidistributing schemes based on
\eqref{eq:sdfdb}, we simply add the terms \eqref{eq:rhoterms2} to the right
hand side of \eqref{eq:sdfdb}. 
For the schemes based on \eqref{eq:sdstabfdb}, on the other hand, for
\eqref{eq:rhoterms} we use the semi-implicit approximation
\begin{equation} \label{eq:rhotermfd}
-  \sum_{p \in \partial_1\mbI} \sliprho^{(p)}\,
(\pol X^m(p)\cdot\pol\ek_1)\,\pol\eta(p)\cdot\pol\ek_2
 - \sum_{p \in \partial_2\mbI} \left\{ [\sliprho^{(p)}]_+ \,
\,(\pol X^{m+1}(p)\cdot\pol\ek_1)
+ [\sliprho^{(p)}]_- \,
\,(\pol X^{m}(p)\cdot\pol\ek_1)\right\}\,\pol\eta(p)\cdot\pol\ek_1,
\end{equation}
where $[ r ]_\pm = \pm \max \{ \pm r, 0 \}$ for $r \in \bR$,
in order to still guarantee unconditional stability. In fact, all of the
adapted schemes still satisfy the theoretical properties which we proved
previously, where for the stability results the surface area $A(\cdot)$ now needs to be replaced by the relevant energy
$E(\cdot)$. For example, we can prove the following theorem, which generalizes
Theorem~\ref{thm:sfvc}.

\begin{thm}[stability and volume conservation]
Let $(\pol X^{m+1},\kappa^{m+1})$ be a solution to \eqref{eq:sdstabfd},
with \eqref{eq:rhotermfd} added to the right hand side of \eqref{eq:sdstabfdb}. 
Then it holds that
\begin{align} \label{eq:stabe2}
E(\pol X^{m+1}) + 2\pi\,\ttau_m\Bigl\langle\pol X^m\cdot\pol\ek_1\,\varkappa_{\rho}^{m+1},~\varkappa_{\rho}^{m+1}|\pol X^m_{\rho}|^{-1}\Bigr\rangle \leq E(\pol X^m).
\end{align}
Moreover, it holds that $M(\pol X^{m+1})=M(\pol X^m)$.
\end{thm}
\begin{proof}
Choosing $\chi = \varkappa^{m+1}$ in \eqref{eq:sdstabfda} and $\pol\eta = \pol X^{m+1} - \pol X^m$ in the adapted \eqref{eq:sdstabfdb} yields
\begin{align}\label{eq:eime2}
&\ttau_m\Bigl\langle\pol X^m\cdot\pol\ek_1,~\varkappa_{\rho}^{m+1},~\varkappa_{\rho}^{m+1}|\pol X^m_{\rho}|^{-1}\Bigr\rangle+\left\langle (\pol X^{m+1}-\pol X^m) \cdot\pol\ek_1, ~|\pol X^{m+1}_\rho|\right\rangle
+ \left\langle (\pol X^m\cdot\pol\ek_1)\,
\pol X^{m+1}_\rho,~(\pol X^{m+1} - \pol X^m)_\rho\, |\pol X^m_\rho|^{-1} \right\rangle \nn\\
&\hspace{0.1em}+\sum_{p \in \partial_1\mbI} \sliprho^{(p)}\,
(\pol X^m(p)\cdot\pol\ek_1)\,(\pol X^{m+1}(p) - \pol X^m(p))\cdot\pol\ek_2
\nn\\
&\hspace{0.1em}
+\sum_{p \in \partial_2\mbI} \left\{\left([\sliprho^{(p)}]_+ \,
\,\pol X^{m+1}(p)
+ [\sliprho^{(p)}]_- \,
\,\pol X^{m}(p)\right)
\cdot\pol\ek_1\right\}\,(\pol X^{m+1}(p) - \pol X^m(p))\cdot\pol\ek_1=0.
\end{align}
Following the proof of Theorem~4.2 in \cite{Barrett19}, we observe that
since $\delta \pol X^{m+1}\in\Vpartial^h$, we have $(\pol X^{m+1}(p)-\pol X^m(p))\cdot\pol\ek_1=0$ for $p\in\partial_1\mbI$, and thus
\begin{align}
 \sum_{p \in \partial_1\mbI} \sliprho^{(p)}\,
(\pol X^m(p)\cdot\pol\ek_1)\,(\pol X^{m+1}(p) - \pol X^m(p))\cdot\pol\ek_2&=\sum_{p\in\partial_1\mbI}\sliprho^{(p)}(\pol X^{m+1}(p)\cdot\pol\ek_1)(\pol X^{m+1}(p)\cdot\pol\ek_2) - \sum_{p\in\partial_1\mbI}(\pol X^m(p)\cdot\pol\ek_1)(\pol X^m(p)\cdot\pol\ek_2)\nn\\
&=\frac{1}{2\pi}\left(A_{\partial_1}(\pol X^{m+1})-A_{\partial_1}(\pol X^m)\right).
\label{eq:sfstableimed2}
\end{align}
In addition, by noting $2\,\beta\,(\beta - \alpha) \geq \beta^2 - \alpha^2$, we get
\begin{align}
&\sum_{p \in \partial_2\mbI} \left\{  [\sliprho^{(p)}]_+ \,
\,\pol X^{m+1}(p)\cdot\pol\ek_1
+ [\sliprho^{(p)}]_- \,
\,\pol X^{m}(p)\cdot\pol\ek_1\right\}\,(\pol X^{m+1}(p) - \pol X^m(p))\cdot\pol\ek_1\nn\\
&\hspace{1cm}\geq\tfrac12 \sum_{p \in \partial_2\mbI} [\sliprho^{(p)}]_+\,
(\pol X^{m+1}(p)\cdot\pol\ek_1)^2\,
- \tfrac12 \sum_{p \in \partial_2\mbI} [\sliprho^{(p)}]_+\,
(\pol X^{m}(p)\cdot\pol\ek_1)^2\nn\\
&\hspace{1.2cm}+  \tfrac12 \sum_{p \in \partial_2\mbI} [\sliprho^{(p)}]_-\,
(\pol X^{m+1}(p)\cdot\pol\ek_1)^2\,
- \tfrac12 \sum_{p \in \partial_2\mbI} [\sliprho^{(p)}]_-\,
(\pol X^{m}(p)\cdot\pol\ek_1)^2\nn\\
&\hspace{1cm}=\frac{1}{2\pi}\left(A_{\partial_2}(\pol X^{m+1}) - A_{\partial_2}(\pol X^m)\right).
\label{eq:sfstableimed3}
\end{align}
Combining \eqref{eq:eime2}, \eqref{eq:sfstableimed1}, \eqref{eq:sfstableimed2}
and \eqref{eq:sfstableimed3} yields the desired stability result 
\eqref{eq:stabe2}. Finally, the volume conservation property can be shown by
choosing $\chi = 1$ in \eqref{eq:sdstabfda} as before.
\end{proof}

\setcounter{equation}{0}
\section{Numerical results}
\label{sec:nr}
In this section we present several numerical experiments to test the fully discrete schemes that we introduced in Section~\ref{sec:fd}. We always employ uniform time steps, $\ttau_m = \ttau$, $m=0,\ldots,M-1$.
 To numerically verify the properties of the discretized schemes, we define the normalized energy function $E(t)/E(0)$, the relative volume loss function $\Delta M(t)$ and the mesh ratio indicator function $\Psi(t)$ via
 \begin{align}
 \left.E(t)/E(0)\right|_{t=t_m}=\frac{E(\pol X^m)}{E(\pol X^0)},\qquad\left.\Delta M(t)\right|_{t=t_m}= \frac{\widetilde M(\pol X^m)-\widetilde M(\pol X^0)}{\widetilde M(\pol X^0)},\nn\\ \left. \Psi(t)\right|_{t=t_m}=\frac{\max_{1\leq j\leq J}|\pol X^m(q_j) - \pol X^m(q_{j-1})|}{\min_{1\leq j\leq J}|\pol X^m(q_j) - \pol X^m(q_{j-1})|},\qquad m\geq 0,\nn
 \end{align}
where $E(\pol X^m)$ is defined in \eqref{eq:E} and
\begin{equation*} 
\widetilde M(\pol X^m)=\pi\left\langle (\pol X^m \cdot\pol\ek_1)^2\, 
\pol\nu^m, \pol\ek_1\,|\pol X^m_\rho|\right\rangle
+ \pi\sum_{p\in\partial_1\mbI} (-1)^{p+1} \,(\pol X^m(p)\cdot\pol\ek_1)^2\,
\pol X^m(p)\cdot\pol\ek_2\,.
\end{equation*}
A proof for $\widetilde M(\pol X^m) = M(\pol X^m)$, in the case 
$\partial\mbI=\partial_0\mbI$, can be found in \cite[(3.10)]{Barrett19}.
In the case of an axisymmetric surface with boundary, the quantity
$\widetilde M(\pol X^m)$, up to an additive constant, represents the volume of
a suitably defined ``interior'' that is bounded by the surface and parts of the
external walls and planes. 

Finally, for computing solutions to the nonlinear systems arising at each time level, we choose the tolerance ${\rm tol}=10^{-10}$ in the Newton's method, recall \eqref{eq:stop}. Then in all our simulations we observe that at most 3 iterations are required in the Newton's method.

\subsection{Numerical results for surface diffusion flow }
\label{sec:sdnr}

\begin{figure}[t]
\centering
\includegraphics[width=0.8\textwidth]{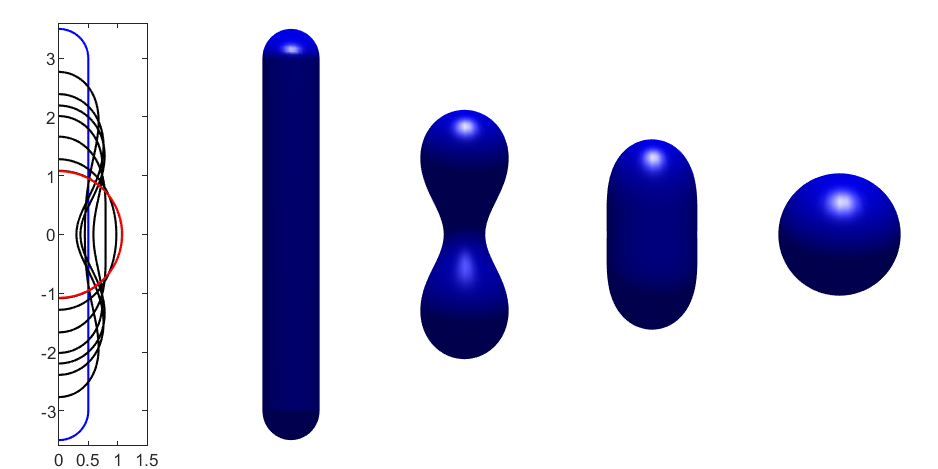}
\caption{Evolution under surface diffusion for a rounded cylinder of dimension $1\times 7\times 1$. The evolution of the generating curves is shown on the left. On the right are the visualizations of the axisymmetric surfaces $\mathcal{S}^m$ at times $t=0, 0.3, 0.5, 1$, respectively. Here $\partial_0\mbI=\{0,1\},~\Delta t=10^{-4}$, $J=128$.}
\label{fig:171E}
\end{figure}

\begin{figure}[!tb]
\centering
\includegraphics[width=0.8\textwidth]{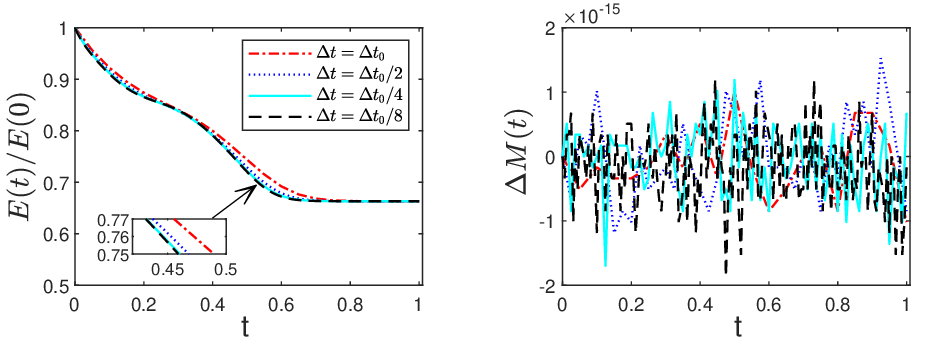}
\caption{The time history of the normalized energy $E(t)/E(0)$ (left panel) and the relative volume loss (right panel) for the evolution under surface diffusion of a rounded cylinder of dimension $1\times 7\times 1$ by using different time step sizes $\Delta t$. Here $\partial_0\mbI=\{0,1\}, ~\Delta t_0 = 0.05$ and $J=128$.}
\label{fig:171Q}
\end{figure}

\vspace{0.5em}
\noindent {\bf Example 1.} We focus on the case of a genus-$0$ surface, so that the boundaries are set to $\partial_0\mbI = \partial\mbI = \{0,1\}$. We first consider the experiment for a rounded cylinder of total dimension $1\times7\times1$. The numerical solutions are computed via the discretized scheme \eqref{eq:sdstabfd}, and the discretization parameters are $\Delta t = 10^{-4}, J=128$. As can be seen in Fig.~\ref{fig:171E},  the cylinder finally reaches a sphere as the stationary solution. Besides, in Fig.~\ref{fig:171Q} we show the time evolution of the normalized energy $E(t)/E(0)$ and the relative volume loss $\Delta M(t)$ by using different time step sizes. We observe the numerical convergence in time for the scheme, and the enclosed volume for the numerical solution is conserved up to the machine precision, confirming our theoretical results.

\begin{figure}[!tbh]
\centering
\includegraphics[width=0.35\textwidth]{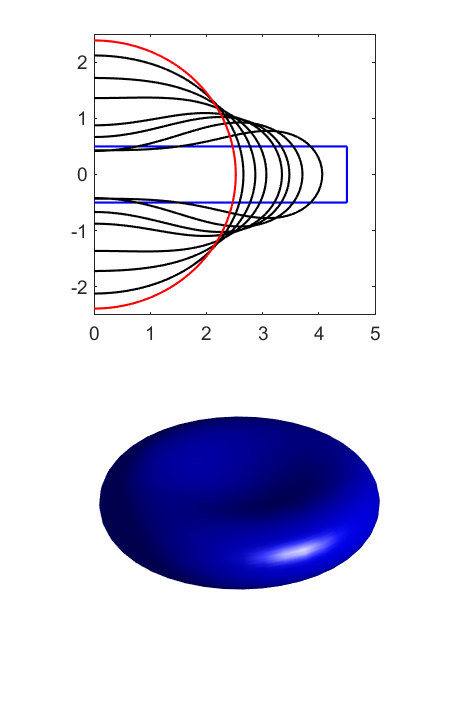}
\includegraphics[width=0.35\textwidth]{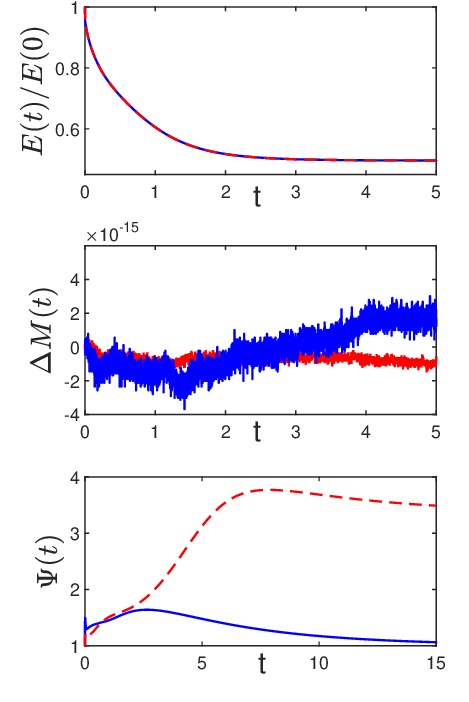}
\caption{Evolution under surface diffusion for a disc of dimension $9\times 1\times 9$ towards the stationary solution. On the left we plot the generating curves at $t=0, 0.2,0.5,0.8,1,1.5,2,3,5$ with the visualizations of the axisymmetric surfaces $\mathcal{S}^m$ at time $t=0.5$. On the right are the plots of the normalized energy $E(t)/E(0)$, the relative volume loss $\Delta M(t)$, and the mesh ratio indicator function $\Psi(t)$ by using \eqref{eq:sdfd} (solid blue line) and \eqref{eq:sdstabfd} (dashed red line). Here $\partial_0\mbI=\{0,1\},~ \Delta t = 10^{-3}$, $J=100$.}
\label{fig:discE}
\end{figure}

In Fig.~\ref{fig:discE}, we consider the experiment for a disc of total dimension $9\times 1\times 9$ via the discretized scheme \eqref{eq:sdfd}, and the discretization parameters are $\Delta t = 10^{-3}, J=128$. The time history of some relevant quantities are plotted, and as a comparison, the plots for numerical solutions computed via scheme \eqref{eq:sdstabfd} are presented together.  We observe the energy dissipation and exact volume conservation for both schemes. For the mesh quality, it can be seen that the mesh ratio indicator tends to 1 for the scheme \eqref{eq:sdfd}, due to the equidistribution property. For the scheme \eqref{eq:sdstabfd}, on the other hand, the mesh ratio indicator reaches larger values but overall does not exceed 4. This indicates that the mesh quality for scheme \eqref{eq:sdstabfd} is preserved well during the simulation although without the equidistribution property.

\begin{figure}[!ht]
\centering
\includegraphics[width=0.35\textwidth]{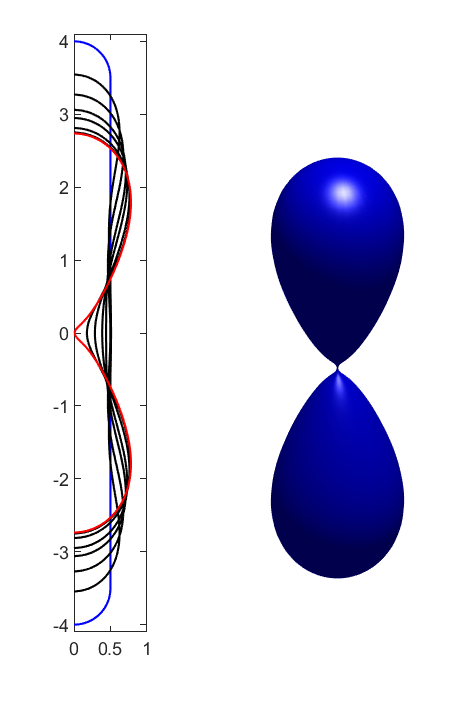}
\includegraphics[width=0.35\textwidth]{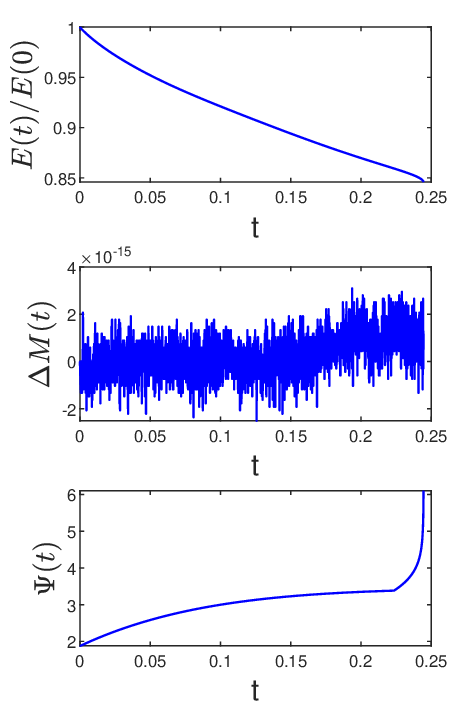}
\caption{Evolution under surface diffusion for a rounded cylinder of dimension $1\times8\times1$ until its pinch off. On the left the generating curves $\Gamma^m$ are shown at times $t=0,0.05,0.1,0.15,0.18,0.22,0.24,0.2452$.
We also visualize the axisymmetric surface $\mathcal{S}^m$ generated by $\Gamma^m$ at time $t=0.2452$. On the right are the plots of the normalized energy, the relative volume loss, and the mesh ratio indicator function. Here $\partial_0\mbI=\{0,1\}, ~J=256, ~\Delta t = 10^{-4}$.}
\label{fig:181E}
\end{figure}

As the numerical results for the schemes \eqref{eq:sdstabfd} and \eqref{eq:sdfd} are often graphically indistinguishable, we only 
visualize the numerical results for the former in the following. We next increase the aspect ratio of the rounded cylinder in Fig.~\ref{fig:171E}, and as initial data we choose a rounded cylinder of total dimension $1\times8\times1$. The discretization parameters are $\ttau = 10^{-4}, J=256$, and the numerical results are visualised in Fig.~\ref{fig:181E}, where we observe that the pinch off event happens at time $t=0.2452$. This is highly consistent with the previous results in \cite{Barrett19, Barrett08JCP}. During the simulations, we also observe the energy dissipation and exact volume conservation.
 
\begin{figure}[t]
\centering
\includegraphics[width=0.35\textwidth]{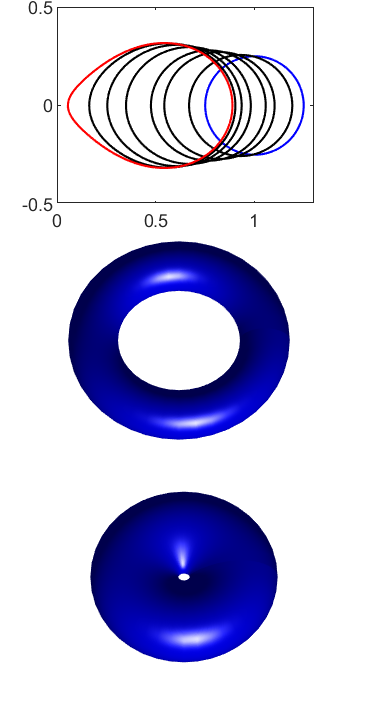}
\includegraphics[width=0.42\textwidth]{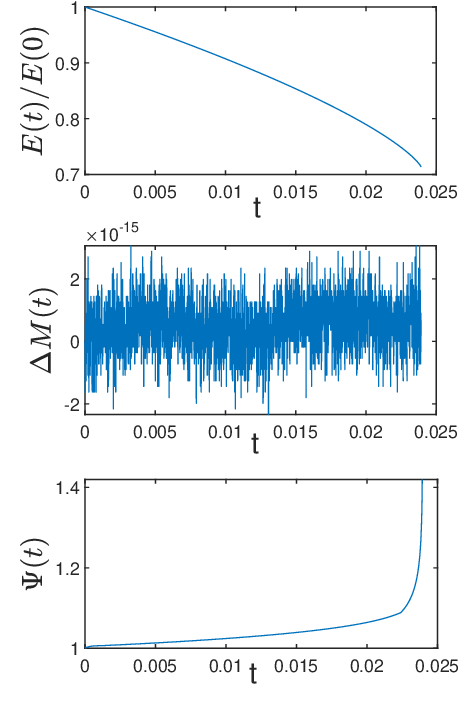}
\caption{Evolution under surface diffusion for a torus with major radius $R=1$ and minor radius $r=0.25$. On the left we show the generating curves $\Gamma^m$ at time $t=0,0.004,0.01,0.013,0.018,0.021,0.023,0.02391$ and the visualizations of the axisymmetric surfaces at $t=0.005, t=0.02391$. On the right are the plots of  the normalized energy $E(t)/E(0)$, the relative volume loss $\Delta M(t)$, and the mesh ratio indicator function $\Psi(t)$. Here $\partial\mbI=\emptyset, ~\Delta t = 10^{-5}$, $J=256$.}
\label{fig:torusE}
\end{figure}

\vspace{0.5em}
\noindent {\bf Example 2.} In this example, we consider the evolution of a torus, so that $\mbI = \RZ$, $\partial\mbI=\emptyset$. As initial data we choose a torus with major radius $R=1$ and minor radius $r =0.25$, and the discretization parameters are $\ttau=10^{-5}, J =256$. The simulation results are illustrated in Fig.~\ref{fig:torusE}, where we observe that the torus shrinks towards the centre and tries to form a genus-$0$ surface, see \cite{Barrett08JCP, Barrett19, Jiang19b}. During the simulations, we observe that the energy dissipation and volume conservation are satisfied well for the numerical solutions. Besides, the mesh ratio indicator remains at small values and increases only slightly when the vertices on $\Gamma^m$ are approaching the $x_2$--axis. 

\begin{figure}[!tbh]
\centering
\includegraphics[width=0.37\textwidth]{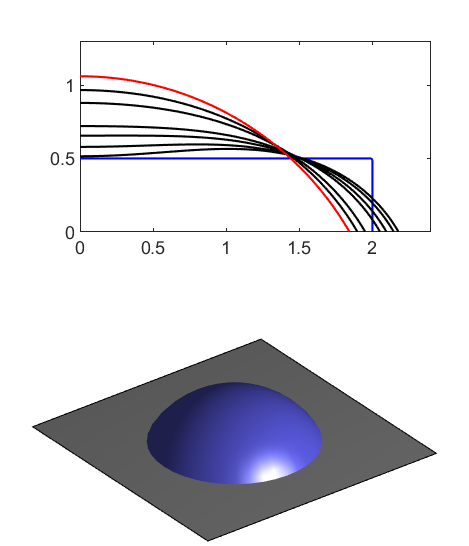}
\includegraphics[width=0.37\textwidth]{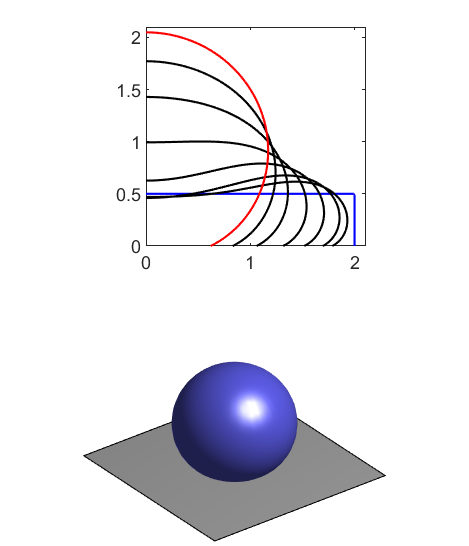}
\caption{Evolution under surface diffusion for a disc droplet of dimension $2\times 1\times 2$ attached to $\mathbb{R}\times\{0\}\times\mathbb{R}$. The generating curves at several times are plotted with $\sliprho^{(1)}=-0.5$ (left panel) and $\sliprho^{(1)}=0.9$ (right panel). We also visualize the axisymmetric surfaces at time $t=1$ on the bottom. Here $\partial_0\mbI=\{0\}, ~\partial_2\mbI=\{1\}, ~\Delta t = 10^{-3}, ~J=128$.}
\label{fig:substrate}
\end{figure}

\begin{figure}[!htb]
\centering
\includegraphics[width=0.80\textwidth]{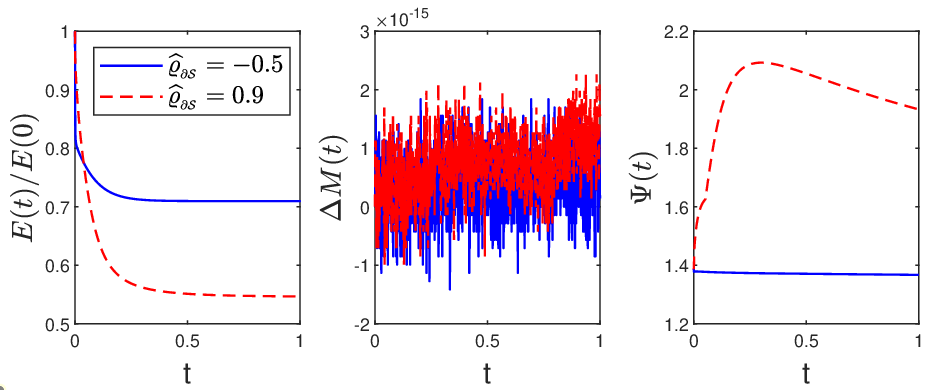}
\caption{The time history of the normalized energy $E(t)/E(0)$, the relative volume loss $\Delta M(t)$, and the mesh ratio indicator function $\Psi(t)$ for the experiment in Fig.~\ref{fig:substrate}.}
\label{fig:subQ}
\end{figure}

\vspace{0.5em}
\noindent{\bf Example 3}. In this example, we are focused on the case when the films are attached to the substrate $\mathbb{R}\times\{0\}\times\mathbb{R}$. We first test the evolution of a small droplet, so that the boundaries are set to $\partial_0\mbI = \{0\}$, $\partial_2 \mbI = \{1\}$. As initial data we choose a disc of total dimension $2\times 1\times 2$, and the discretization parameters are $\ttau=10^{-3}, J=128$. The simulation results for $\sliprho^{(1)}=-0.5$ and $0.5$  are reported in Fig.~\ref{fig:substrate}. When $\sliprho^{(1)}=-0.5$, we observe that the droplet eventually forms a spherical shape with an acute contact angle. On the other hand, when $\sliprho^{(1)}=0.9$, the droplet forms a spherical shape with an obtuse contact angle. During the simulations, the energy dissipation and the volume conservation are observed for the numerical solutions. The mesh ratio indicator remains at small values for both cases. 

\begin{figure}[!htb]
\centering
\includegraphics[width=0.85\textwidth]{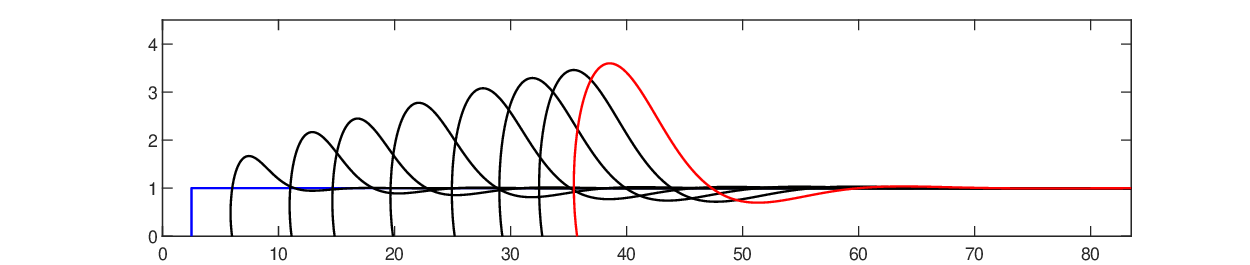}
\includegraphics[width=0.75\textwidth]{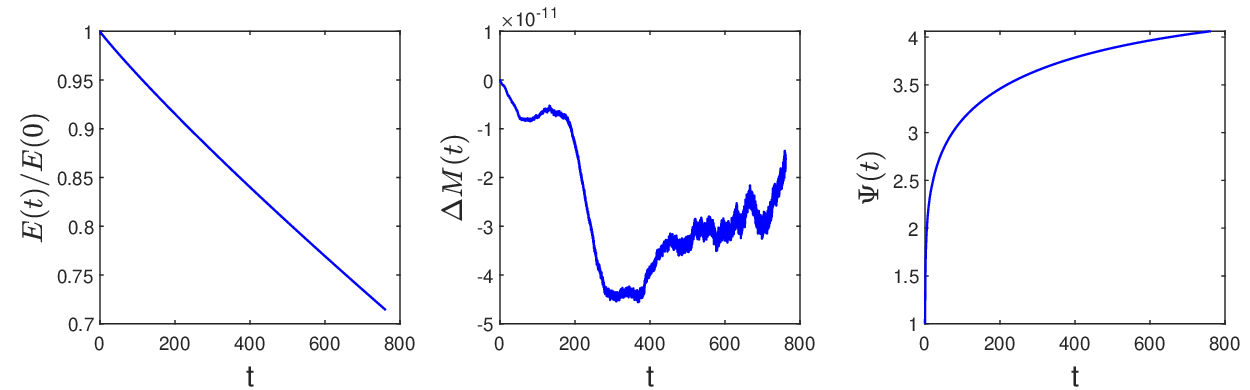}
\caption{Evolution under surface diffusion for a large disc of dimension $83.5\times1\times 83.5$ with a hole of dimension $2.5\times1\times 2.5$. On the top are plots of the generating curves $\Gamma^m$ at times $t=0,10,50,100,200,350, 500, 650,800$. On the bottom are  plots of the normalized energy, the relative volume loss and the mesh ratio indicator function. Here $\partial_2\mbI=\{0\},~\partial_{_D}\mbI=\{1\}, ~\sliprho^{(0)}=-0.5, ~J = 820,~ \ttau=10^{-3}$.}
\label{fig:holeg}
\end{figure}

We next study the evolution of a large disc with a relatively small hole. As initial data we choose a disc of dimension $83.5\times 1\times 83.5$ with a hole of dimension $2.5\times1\times 2.5$ in the centre. For the boundaries we set $\partial_2\mbI=\{0\}$, $\partial_{_D}\mbI =\{1\}$ with $\sliprho^{(0)}=-0.5$  to simulate the hole growth in solid-state dewetting \cite{Wong00, Wang15,Jiang19a}. The discretization parameters are $J=820, \ttau=10^{-3}$, and the numerical results are shown in Fig.~\ref{fig:holeg}. In this experiment, we observe that the hole gets larger and the retracting edge forms a
thickened ridge followed by a valley. During the evolution, it can be seen that the ridge gradually rises and the valley sinks. This is consistent to the results in \cite{Wong00, Wang15}. Besides, we observe that the energy dissipation and the volume conservation are well satisfied for the numerical solutions, and the mesh quality is well preserved. 

\begin{figure}[!ht]
\centering
\includegraphics[width=0.75\textwidth]{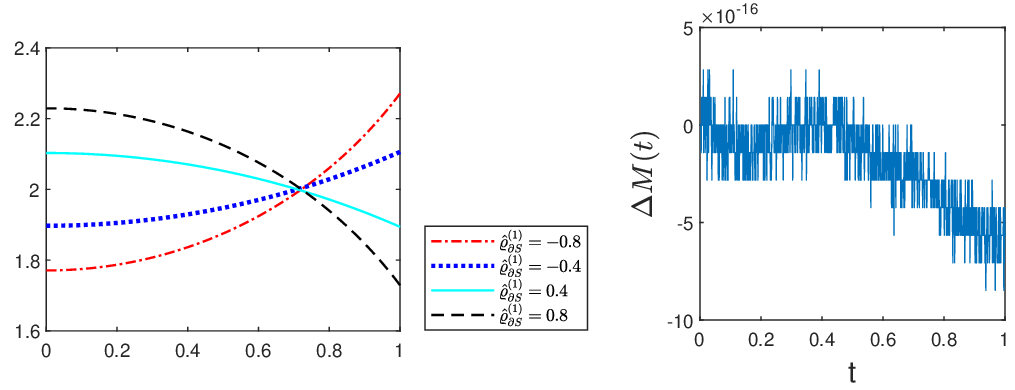}
\caption{Evolution under surface diffusion for a disc attached to an infinite cylinder of radius 1. On the left are the generating curves of the stationary solutions for different $\sliprho^{(1)}$ and on the right we plot the time history of the relative volume loss for $\sliprho^{(1)}=-0.8$. Here $\partial_0\mbI=\{0\},~ \partial_1\mbI=\{1\}, ~J=128,~ \ttau=10^{-3}$.}
\label{fig:icylinder}
\end{figure}

\begin{figure}[!ht]
\centering
\includegraphics[width=0.75\textwidth]{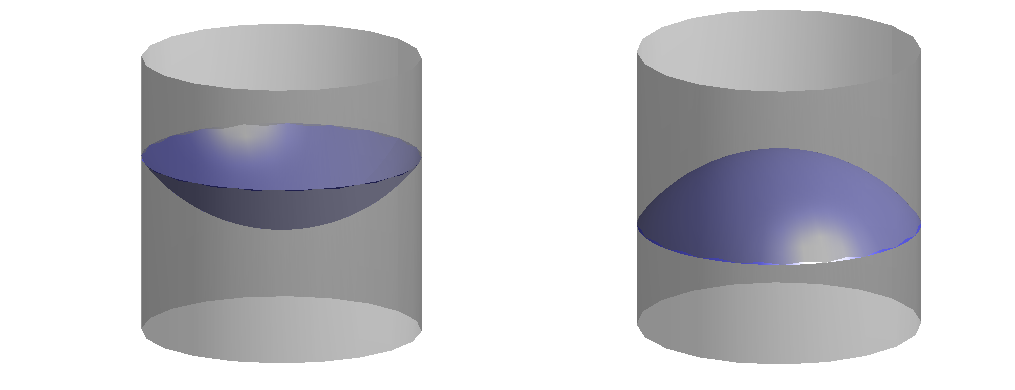}
\caption{Visualisations of the axisymmetric surfaces for $\sliprho^{(1)}=-0.8$ (left panel) and $\sliprho^{(1)}=0.8$ (right panel) for the experiment in Fig.~\ref{fig:icylinder}. }
\label{fig:icylinderEQ}
\end{figure}

\vspace{0.5em}
\noindent{\bf Example 4}. For the next experiment, we consider a disc attached to an infinite cylinder of radius 1, with prescribed contact angle conditions. The boundaries are set to $\partial_0\mbI=\{0\},~\partial_1\mbI=\{1\}$. We start with a disc of total dimension $2\times1\times 2$, and the discretization parameters are $J=128, \ttau = 10^{-3}$. The generating curves of the stationary solutions are presented in Fig.~\ref{fig:icylinder} under four different parameters $\sliprho^{(1)}=-0.8,-0.4,0.4,0.8$. We also plot the time history of the relative volume loss  for the case $\sliprho^{(1)}=-0.8$, and observe the exact volume conservation for the numerical solutions. For this experiment, the corresponding axisymmetric surfaces for $\sliprho^{(1)}=-0.8,\sliprho^{(1)}=0.8$ are visualized in Fig.~\ref{fig:icylinderEQ}.

\begin{figure}[t]
\centering
\includegraphics[width=0.85\textwidth]{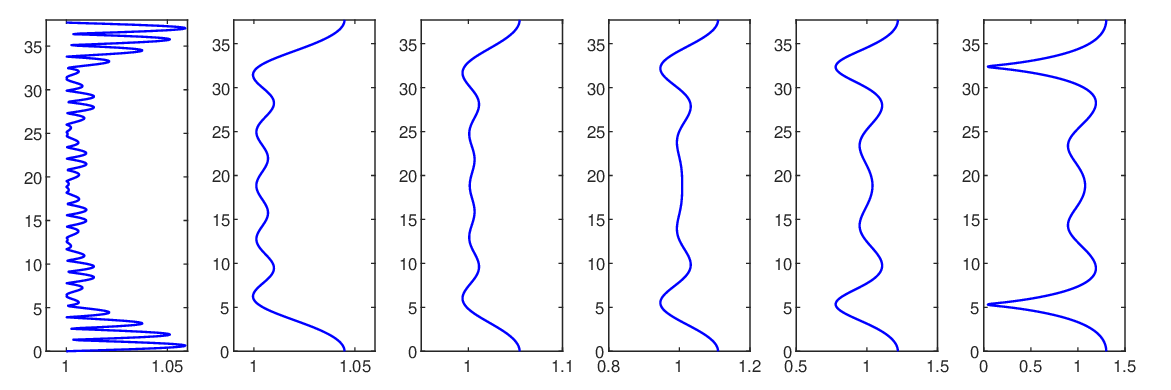}
\includegraphics[width=0.85\textwidth]{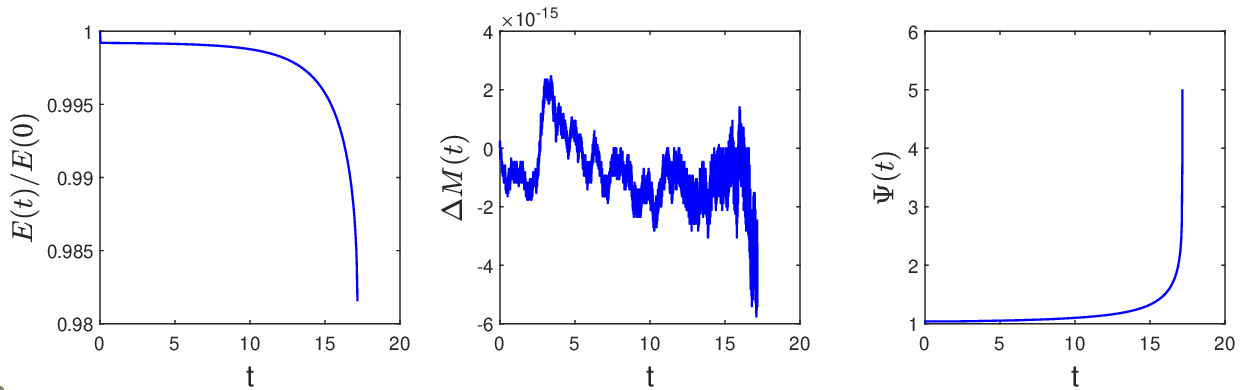}
\caption{Evolution under surface diffusion for a cylindrical shape attached to two hyperplanes under small perturbations. On the top are the plots of the generating curves at times $t=0, 2, 4,10,15, 17.1543$. On the bottom are the plots of the normalized energy, the relative volume loss and the mesh ratio indicator during the simulations. Here $\partial\mbI=\partial_2\mbI=\{0,1\},~\sliprho^{(0)}=\sliprho^{(1)} =0, ~J=512, ~\ttau = 10^{-4}$.}
\label{fig:stability}
\end{figure}

\vspace{0.5em}
\noindent{\bf Example 5}. We end this subsection by testing the stability of a cylinder under small perturbations \cite{ColemanFM96}. As initial data we choose 
\begin{align}
r(t) = 1 + 0.01\times|\sin(2\,z) + \sin(13z/6) + \sin(7z/3) + \sin(5z/2) + \sin(8z/3)+ \sin(17z/6)|, \quad 0\leq z\leq 12\pi,\nn
\end{align}
which forces a small amplitude perturbation of the cylinder. The cylindrical shapes are attached to two parallel hyperplanes $\mathbb{R}\times\{0\}\times\mathbb{R}$ and $\mathbb{R}\times\{12\pi\}\times\mathbb{R}$, so that we have 
$\partial_2\mbI=\{0,1\}.$ We set $\sliprho^{(0)}=\sliprho^{(1)} =0$, and the discretization parameters are $J=512, \ttau=10^{-4}$. The numerical results are reported in Fig.~\ref{fig:stability}, where we observe the final breakup of the cylinder at time $t=17.1543$ although the wavelengths of the initial perturbations are small. This is consistent to the results in \cite{ColemanFM96}. During the simulation, the energy dissipation and volume conservation are observed to be well satisfied, and the mesh ratio indicator remains at small values during the evolution except when the vertices of $\Gamma^m$ are approaching the $x_2$-axis.

\subsection{Numerical results for the conserved mean curvature flow}
\begin{figure}[t]
\centering
\includegraphics[width=0.35\textwidth]{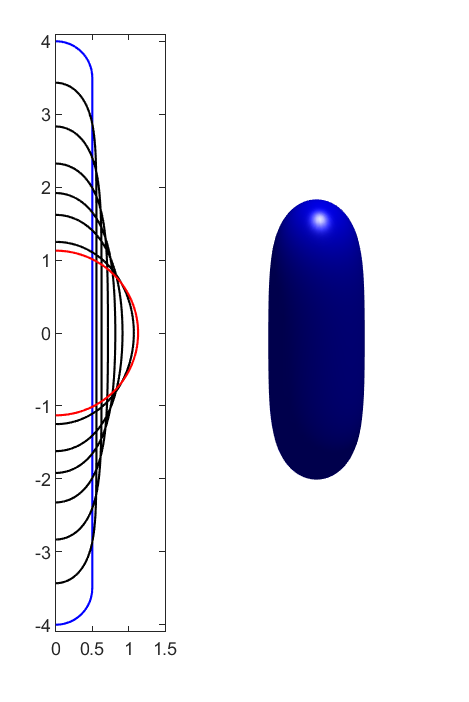}
\includegraphics[width=0.35\textwidth]{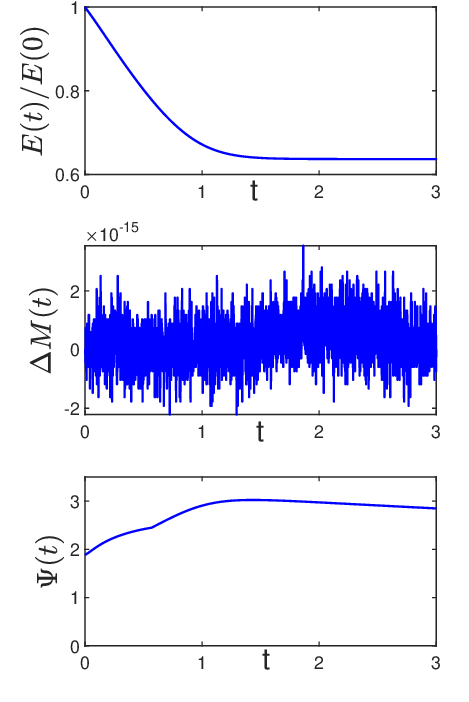}
\caption{Evolution under conserved mean curvature flow for a rounded cylinder of dimension $1\times8\times1$ towards the stationary solution. On the left are the plots of the generating curves $\Gamma^m$ at times $t=0,0.2,0.4,0.6,0.8,1,1.5,4$.
We also visualize the axisymmetric surface $\mathcal{S}^m$ generated by $\Gamma^m$ at time $t=0.6$. On the right are the plots of the normalized energy, the relative volume loss, and the mesh ratio indicator function. Here $\partial\mbI=\partial_0\mbI=\{0,1\}, ~J=256, ~\Delta t= 10^{-3}$.}
\label{fig:MCF181}
\end{figure}

\vspace{0.5em}
\noindent{\bf Example 6}.  For the conserved mean curvature flow,  we consider two numerical experiments that have been investigated for surface diffusion flow in \S\ref{sec:sdnr}. We start by testing the evolution of a cylinder with the same numerical setup as in Fig.~\ref{fig:181E}, except that we choose $\ttau=10^{-3}$. The numerical results are reported in Fig.~\ref{fig:MCF181}, where we observe that the cylinder tries to form a sphere as the stationary solution, while in Fig.~\ref{fig:181E}, the formation of a singularity is observed for surface diffusion. In particular, we find that the generating curves remain convex during the simulation, and thus do the axisymmetric surfaces. This is in line with the theoretical result from \cite{Huisken1987volume}, where it is shown that convexity is preserved for conserved mean curvature flow.
This is in contrast to surface diffusion and the intermediate evolution law,
for which convexity in general is not preserved, see 
\cite{GigaI99,Ito02,EscherI05}, as well as Figs.~\ref{fig:171E}, 
\ref{fig:im181xi1} and \ref{fig:im181xi6}.
The exact volume conservation and good mesh quality are also observed for the numerical solutions in Fig.~\ref{fig:MCF181}.

\begin{figure}[!ht]
\centering
\includegraphics[width=0.35\textwidth]{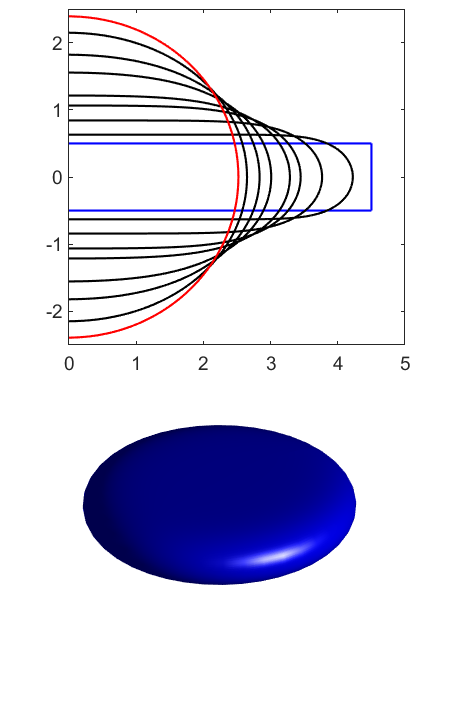}
\includegraphics[width=0.35\textwidth]{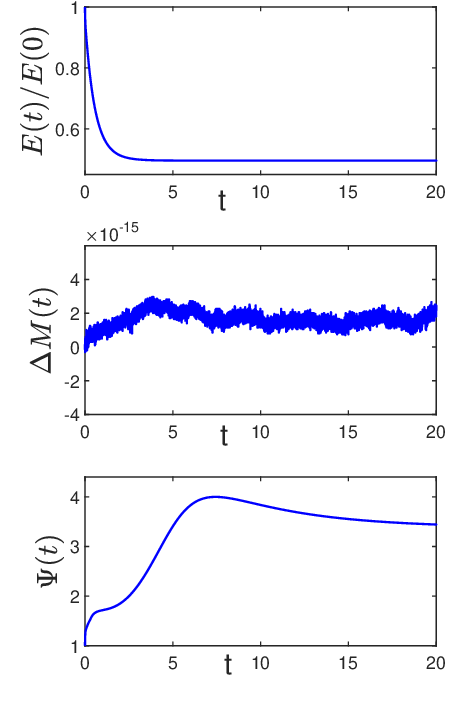}
\caption{Evolution under conserved mean curvature flow for a disc of dimension $9\times 1\times 9$ towards the stationary solution. One the left are the plots of the generating curves $\Gamma^m$ at $t=0, 0.2,0.5,0.8,1,1.5,2,3,20$ together with a visualization of the axisymmetric surface $\mathcal{S}^m$ generated by $\Gamma^m$ at time $t=0.5$. On the right are the plots of the normalized energy, the relative volume loss, and the mesh ratio indicator versus time. Here $\partial\mbI=\partial_0\mbI=\{0,1\},~ \Delta t = 10^{-3}$, $J=100$.}
\label{fig:MCFdiscE}
\end{figure}

Moreover, we consider the evolution of a disc with the same numerical setup as in Fig.~\ref{fig:discE}, and the simulation results are illustrated in Fig.~\ref{fig:MCFdiscE}. In this experiment, the disc remains convex during the simulation and finally forms a sphere as the stationary solution. Besides, we observe that the energy dissipation and volume conservation are satisfied for the numerical solutions, and the mesh quality is generally well preserved.

\subsection{Numerical results for the intermediate evolution flow}
\begin{figure}[t]
\centering
\includegraphics[width=0.45\textwidth]{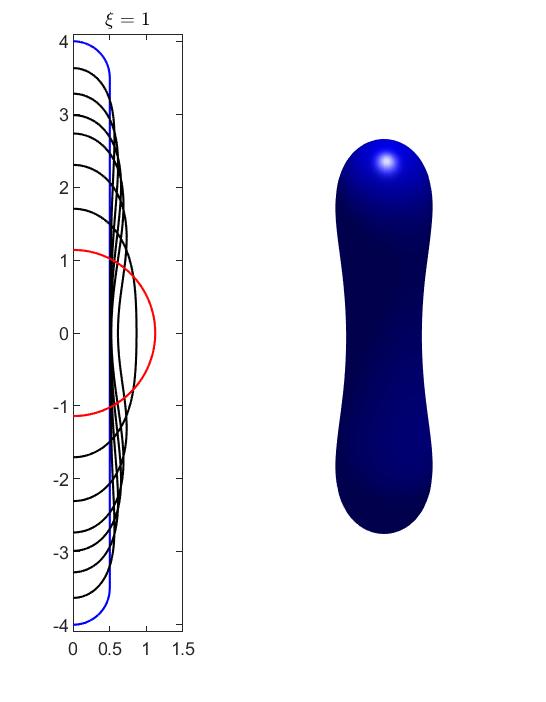}
\includegraphics[width=0.4\textwidth]{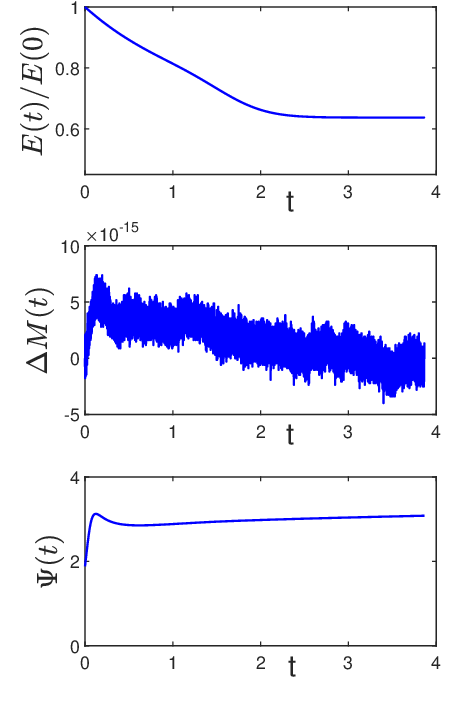}
\caption{Evolution under the intermediate evolution flow for a rounded cylinder of dimension $1\times8\times1$ towards the stationary solution. On the left are plots of the generating curves $\Gamma^m$ at $t=0,0.2,0.4,0.6,0.8,1.2,1.8,3.5$. We also visualize the axisymmetric surface $\mathcal{S}^m$ generated by $\Gamma^m$ at $t=0.6$. On the right are plots of the normalized energy, the relative volume loss, and the mesh ratio indicator function. Here $\partial\mbI=\partial_0\mbI=\{0,1\},~ J=256, ~\Delta t= 10^{-4}$, $\alpha=1$ and $\xi=1$.}
\label{fig:im181xi1}
\end{figure}

\begin{figure}[!htb]
\centering
\includegraphics[width=0.4\textwidth]{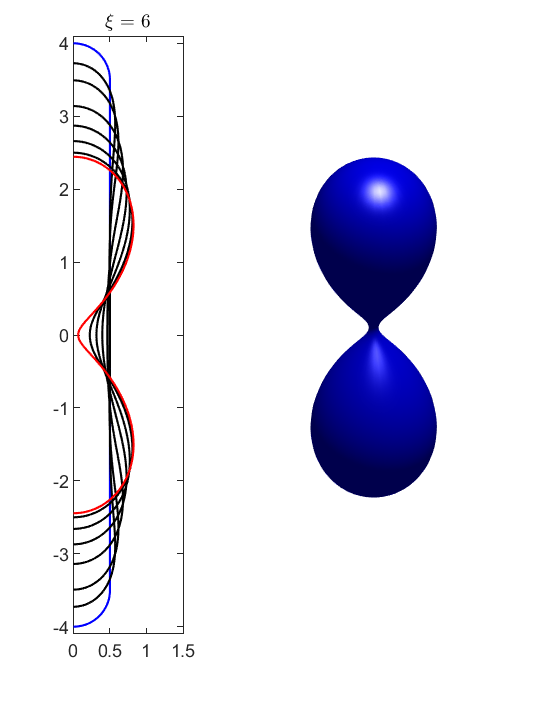}
\includegraphics[width=0.35\textwidth]{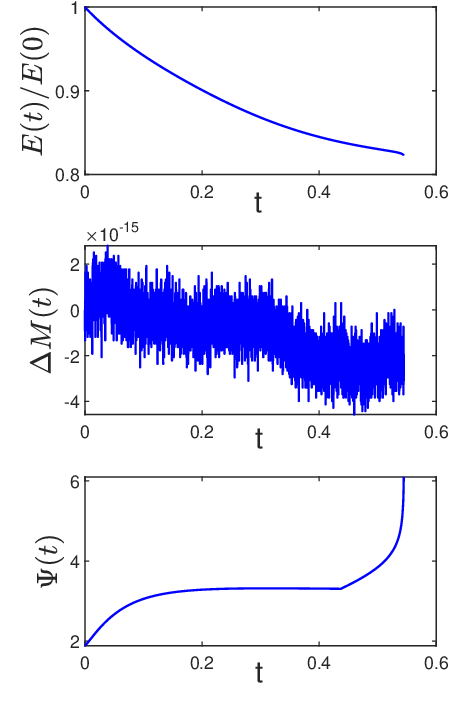}
\caption{Evolution under the intermediate evolution flow for a rounded cylinder of dimension $1\times8\times1$ until its breakup. On the left the generating curves $\Gamma^m$ are shown at $t=0,0.05,0.1,0.2,0.3,0.4,0.5,0.544$.
We also visualize the axisymmetric surface $\mathcal{S}^m$ generated by $\Gamma^m$ at $t=0.544$. On the right are the plots of the normalized energy, the relative volume loss, and the mesh ratio indicator function versus time. Here $\partial\mbI=\partial_0\mbI=\{0,1\}, ~J=256, ~\Delta t= 10^{-4}$, $\alpha=1$ and $\xi=6$.}
\label{fig:im181xi6}
\end{figure}

\vspace{0.5em}
\noindent{\bf Example 7}:
For the intermediate evolution flow, we can repeat all the examples in  \S\ref{sec:sdnr}. For simplicity, here we only show an example with the same numerical setup as in Fig.~\ref{fig:181E}. We first fix $\alpha=1$ and choose $\xi=1$, so that the flow interpolates between surface diffusion and conserved mean curvature flow. As shown in Fig.~\ref{fig:im181xi1}, we observe that the cylinder finally forms a sphere as the stationary solution. We then set $\xi=6$ so that surface diffusion gains more weight in the evolution flow. In this experiment, as shown in Fig.~\ref{fig:im181xi6}, we observe the occurrence of pinch off at $t=0.544$ although the evolution is much slower compared to that in Fig.~\ref{fig:181E}. During the simulations, we observe the energy dissipation and the exact volume conservation, and the mesh quality is generally well preserved except when the cylinder is about to break up.

\begin{figure}[!htb]
\centering
\includegraphics[width=0.13\textwidth]{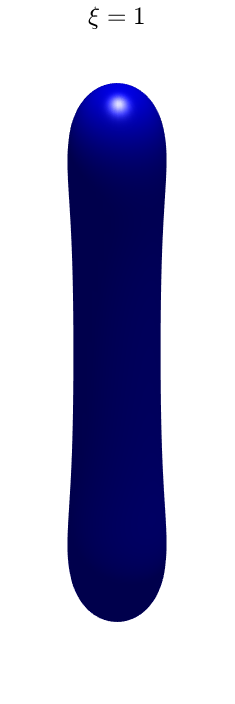}
\includegraphics[width=0.13\textwidth]{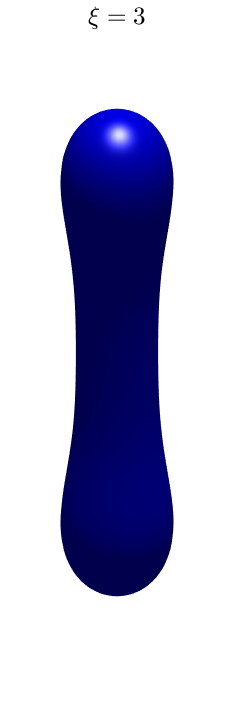}
\includegraphics[width=0.13\textwidth]{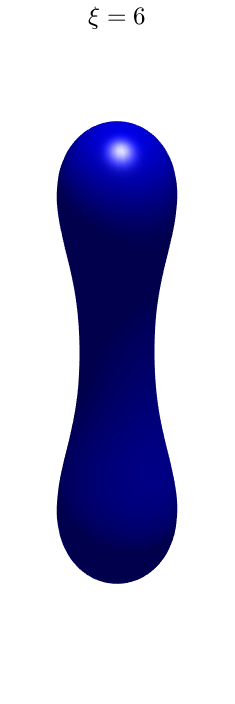}
\includegraphics[width=0.13\textwidth]{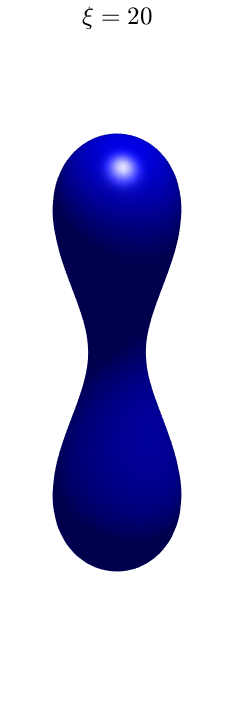}
\includegraphics[width=0.13\textwidth]{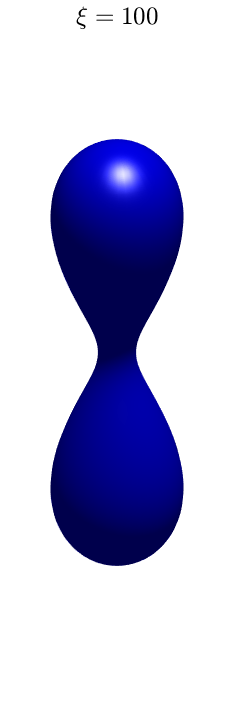}
\includegraphics[width=0.13\textwidth]{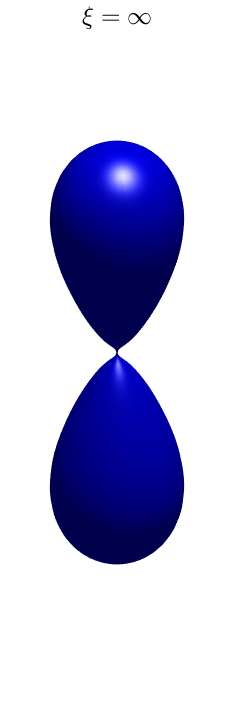}\\
\includegraphics[width=0.13\textwidth]{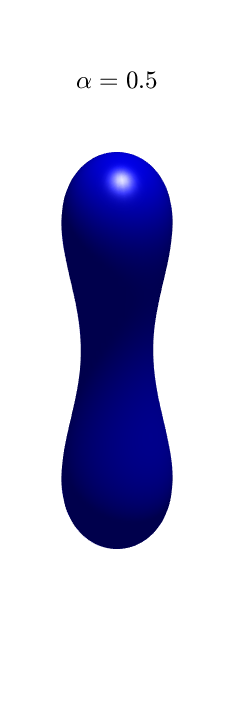}
\includegraphics[width=0.13\textwidth]{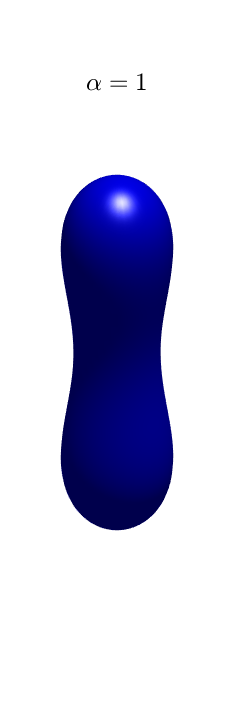}
\includegraphics[width=0.13\textwidth]{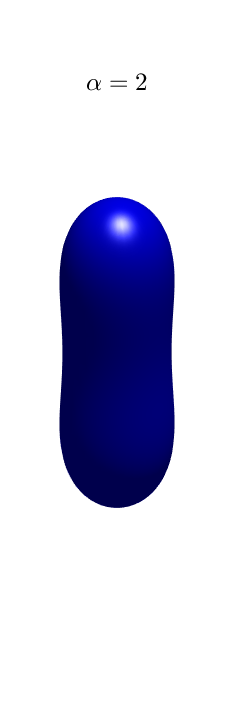}
\includegraphics[width=0.13\textwidth]{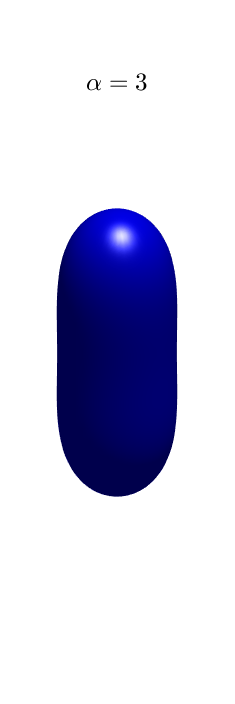}
\includegraphics[width=0.13\textwidth]{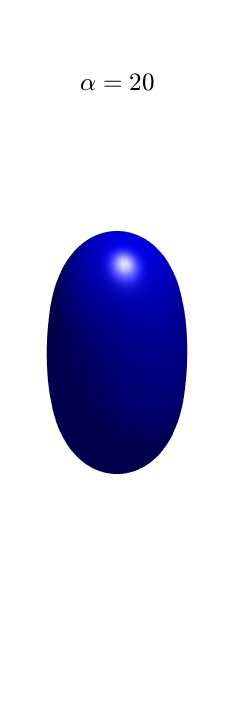}
\includegraphics[width=0.13\textwidth]{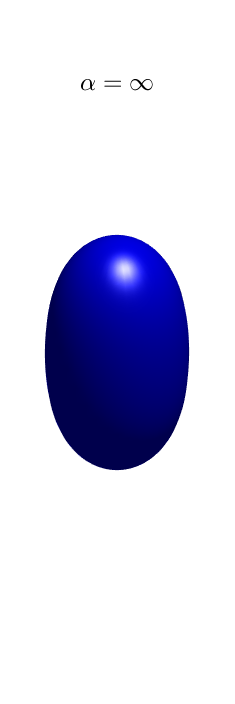}
\caption{Evolution under the intermediate evolution flow for a rounded cylinder of dimension $1\times 8\times 1$ with different pairs of $(\xi,~\alpha)$. Top panel: we fix $\alpha=1$ and visualize the axisymmetric surfaces at $t=0.2452$ for different values of $\xi$. Bottom panel: we fix $\xi=1$ and visualize the axisymmetric surfaces at $t=1.0$ for different values of $\alpha$. Here $\partial\mbI=\partial_0\mbI=\{0,1\}, ~J=256, ~\Delta t= 10^{-4}$. Note that $\alpha=1,\xi=+\infty$ represents the evolution under surface diffusion in Fig.~\ref{fig:181E}, and $\alpha=+\infty,\xi=1$ represents the evolution under conserved mean curvature flow in Fig.~\ref{fig:MCF181}. }
\label{fig:imtm}
\end{figure}

\begin{figure}[!htb]
\centering
\includegraphics[width=0.88\textwidth]{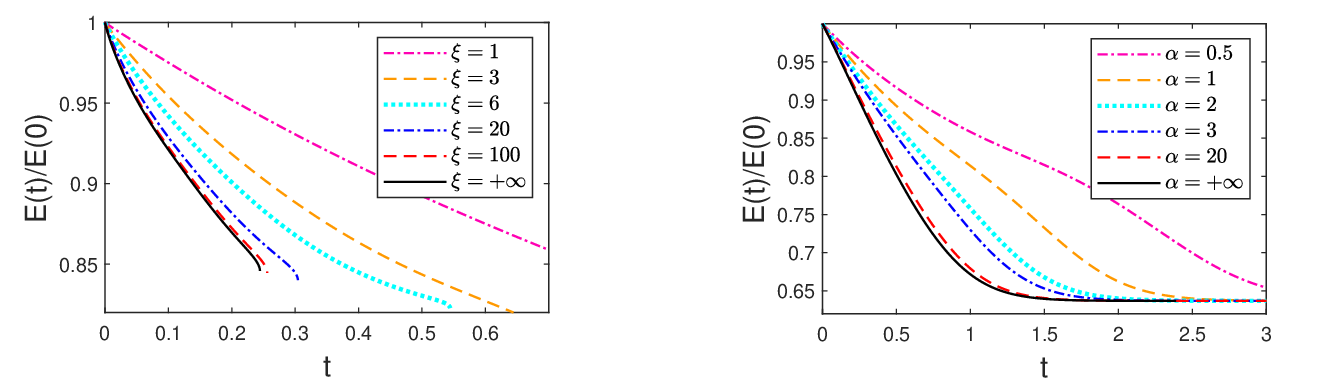}
\caption{The time history of the normalized energy under different pairs of $(\xi,~\alpha)$ for the experiments in Fig.~\ref{fig:imtm}. Left panel: we fix $\alpha=1$ and choose different values of $\xi$. The pinch-off occurs for $\xi=6,20,100,+\infty$ at $t_c=0.5440, 0.3044,0.2559,0.2452$, respectively. Right panel: we fix $\xi=1$ and choose different values of $\alpha$.}
\label{fig:imtm1}
\end{figure}

To further investigate the intermediate evolution flow, we use the same numerical setup as above but consider the evolution when $\xi$ or $\alpha$ approaches infinity. Precisely, we fix $\alpha=1$ and visualize the axisymmetric surfaces at $t=0.2452$ under different values of $\xi$, as illustrated in Fig.~\ref{fig:imtm} (top panel).  Similarly, the axisymmetric surfaces at $t=1.0$ under different values of $\alpha$ are visualised in Fig.~\ref{fig:imtm} (bottom panel). As a comparison, we also include results for the experiments in Fig.~\ref{fig:181E} and Fig.~\ref{fig:MCF181}, which correspond to surface diffusion flow ($\alpha=1,\xi=+\infty$) and conserved mean curvature flow ($\alpha=+\infty,\xi=1$), respectively. As $\xi$ increases, we observe that the axisymmetric surfaces tend to break up more easily. Besides, with higher values of $\alpha$, the axisymmetric surfaces at time $t=1.0$  tend to become convex and mimic the property of convexity preservation that the conserved mean curvature flow satisfies. This confirms the theoretical
results from \cite{EscherGI02}. The corresponding time history of the normalized energy for the experiments is shown in Fig.~\ref{fig:imtm1}, where similar convergences are observed as well.

\section{Conclusion}
\label{sec:con}
In this work, we proposed several volume-preserving parametric finite element schemes for discretizing surface diffusion flow, conserved mean curvature flow and the intermediate evolution flow in an axisymmetric setting. The proposed schemes are based on variational discretizations of weak formulations that allow for tangential degrees of freedoms. A suitable weighting between old and new time level of a discrete normal allowed us to prove exact volume conservation of the numerical solutions. Some of the numerical methods enjoy the desirable property of unconditional stability, while others exhibit an asymptotic equidistribution property. Numerical results were presented to demonstrate the accuracy and efficiency of the proposed schemes, and to numerically verify these good properties. 

\section*{Acknowledgement}
The work of Bao was supported by the Ministry of Education of Singapore grant MOE2019-T2-1-063 (R-146-000-296-112). The work of Zhao was funded by the Alexander von Humboldt Foundation. 

\bibliographystyle{model1b-num-names}
\bibliography{thebib}
\end{document}